\documentclass[a4paper,reqno,draft,12pt]{amsart}
\usepackage{amsmath}
\usepackage{amsbsy}
\usepackage{amsopn}
\usepackage{upref}
\usepackage{amsthm}
\usepackage{amsfonts}
\usepackage{amssymb}
\usepackage{mathrsfs}
\usepackage{times}
\newenvironment{bew}[2]{\removelastskip\vspace{6pt}\noindent
 {\it Proof  #1.}~\rm#2}{\par\vspace{6pt}}  
\allowdisplaybreaks
\numberwithin{equation}{section}
\allowdisplaybreaks
\theoremstyle{plain}
 \newtheorem{theorem}{Theorem}[section]
 
 \newtheorem{proposition}[theorem]{Proposition}
 \newtheorem{lemma}[theorem]{Lemma}
\theoremstyle{definition}
 \newtheorem{definition}{Definition}[section]
\theoremstyle{remark}
 \newtheorem{remark}{Remark}
\newcommand{\ep}{\varepsilon}
\newcommand{\p}{\partial}
\newcommand{\RR}{\mathbb{R}}

\begin{document}
\title
[Dispersive Flow]
{
The initial value problem for a third-order dispersive flow
into compact almost Hermitian manifolds}
\author
[E.~Onodera]
{Eiji ONODERA}
\address{Mathematical Institute, Tohoku University, Sendai 980-8578, Japan}
\email{sa3m09@math.tohoku.ac.jp}
\subjclass[2000]
{Primary 35Q55; Secondary 35Q53, 53C44}
\keywords{Schr\"odinger map, geometric analysis, energy
method, smoothing effect}
\thanks{The author is supported by the JSPS Research Fellowships 
        for Young Scientists and 
        the JSPS Grant-in-Aid for Scientific Research No. 19$\cdot$3304.}
\maketitle
\begin{abstract}
We present a time-local existence theorem 
of the initial value problem for 
a third-order dispersive evolution equation 
for open curves on  compact almost Hermitian manifolds 
arising in the geometric analysis of vortex filaments. 
This equation causes the so-called loss of 
one-derivative since the target manifold is 
not supposed to be a K\"ahler manifold. 
We overcome this difficulty by using a 
gauge transformation of a multiplier
on the pull-back bundle 
to eliminate the bad first order 
terms essentially. 
\end{abstract}

\section{Introduction}
\label{section:introduction}
Let $(N, J, g)$ be a compact almost Hermitian manifold 
with an almost complex structure $J$ and a hermitian metric $g$, 
and let $\nabla$ be the Levi-Civita connection with respect to $g$. 
$X$ denotes $\RR$ or $\RR/\mathbb{Z}$. 
Consider the initial value problem of the form
\begin{alignat}{2}
  u_{t}
& = 
  a\,\nabla_x^2u_x
  + J_{u}\nabla_xu_x
  +b\,g_{u}(u_x, u_x)u_x
&
  \quad\text{in}\quad
& \RR\times X,
\label{equation:pde}
\\
  u(0,x)
& =
  u_0(x)
&
  \quad\text{in}\quad
& X,
\label{equation:data}
\end{alignat}
where 
$a, b\in \RR$ are constants, 
$u(t,x)$ is an $N$-valued unknown function of 
$(t,x)\in \RR\times X$, 
$u_t(t,x)
=du_{(t,x)}( \left(\p/\p{t}\right)_{(t,x)})$, 
$u_x(t,x)
=du_{(t,x)}( \left(\p/\p{x}\right)_{(t,x)})$,  
$du_{(t,x)}:T_{(t,x)}(\mathbb{R}\times X)\to T_{u(t,x)}N$ 
is the differential of the mapping
$u$ at $(t,x)$, 
$\nabla_x$
is the covariant derivative
induced from $\nabla$ with respect to $x$ 
along the mapping $u$,
and 
$J_u$ and $g_u$ mean 
the almost complex structure and the metric at $u{\in}N$ respectively. 
The equation \eqref{equation:pde} is 
an equality of sections of the pull-back bundle $u^{-1}TN$. 
We call the solution of \eqref{equation:pde} 
a dispersive flow. 
In particular, when $a=b=0$, 
this is called a one-dimensional Schr\"odinger map.  
\par
Examples of dispersive flows arise in classical mechanics: 
the motion of vortex filament, 
the Heisenberg ferromagnetic spin chain and etc. 
Solutions to these physical models are valued in 
two-dimensional unit sphere $\mathbb{S}^2\subset \RR^3$.
For $\vec{u}=(u_1,u_2,u_3)\in\RR^3$ and 
$\vec{v}=(v_1,v_2,v_3)\in\RR^3$, 
let 
$$
\vec{u}\cdot\vec{v}=u_1v_1+u_2v_2+u_3v_3, 
\quad 
\lvert\vec{u}\rvert=\sqrt{\vec{u}\cdot\vec{u}}, 
$$
$$
\vec{u}\times\vec{v}
=
(u_2v_3-u_3v_2,u_3v_1-u_1v_3,u_1v_2-u_2v_1). 
$$
In \cite{DR}, Da Rios formulated the equation modeling 
the motion of vortex filament of the form 
\begin{equation}
\vec{u}_t
= 
\vec{u} \times \vec{u}_{xx}, 
\label{equation:darios}
\end{equation}
where $\vec{u}(t,x)\in \mathbb{S}^2$ denotes 
the velocity vector along the
space curve describing the position of 
the vortex filament in $\RR^3$ at $(t,x)$, 
$t$ is the time and $x$ is the arc-length in this physical model. 
See also, e.g.,
\cite{Hasimoto} and \cite{LL} for physical backgrounds of 
\eqref{equation:darios}. 
The physical model \eqref{equation:darios} is 
an example of the equation of the one-dimensional Schr\"odinger map. 
Our equation \eqref{equation:pde} 
with $b=a/2$ 
geometrically generalizes 
an $\mathbb{S}^2$-valued physical model 
\begin{equation}
\vec{u}_t
= 
\vec{u} \times \vec{u}_{xx}
+
a
\left[
\vec{u}_{xxx} 
+
\frac{3}{2}
\{\vec{u}_x\times (\vec{u}\times \vec{u}_x )\}_x
\right]
\label{equation:NT}
\end{equation}
describing the motion of vortex filament in $\RR^3$
proposed by Fukumoto and Miyazaki in \cite{FM}. 
\par
Here we state the known results on the mathematical analysis 
of the IVP \eqref{equation:pde}-\eqref{equation:data}. 
There has been many studies on
the existence of solutions to 
\eqref{equation:pde}-\eqref{equation:data} 
both on $X=\RR$ and $\RR/ \mathbb{Z}$
only when  
$(N, J, g)$ is a K\"ahler manifold. 
See 
\cite{CSU}, 
\cite{Ding}, 
\cite{Koiso}, 
\cite{MCGAHAGAN}, 
\cite{NSU1}, 
\cite{NSU2}, 
\cite{PWW2},
\cite{SSB} 
for $a=0$ 
and 
\cite{NT}, 
\cite{Onodera1}, 
\cite{Onodera2}, 
\cite{TN} 
 for $a\ne0$. 
Time-local existence theorems were proved by 
some classical energy estimates with respect to 
the following quantity like the $L^2$-energy
$$
\left\|V\right\|_{L^2(X;TN)}^2
=
\int_{X}g_{u(x)}
\left(
V(x),V(x)
\right)dx
\quad
\text{for}
\quad
V\in \Gamma(u^{-1}TN). 
$$
More precisely, if 
$\nabla$ is a metric connection ($\nabla{g}=0$) 
and $g$ is a K\"ahler metric ($\nabla{J}=0$), 
then the equation \eqref{equation:pde} behaves like 
symmetric hyperbolic systems, 
and the classical energy method works well. 
This fact is closely related with 
the geometric studies of the good structure of 
the equation of dispersive flow into 
a compact Riemann surface on $\mathbb{R}$. 
Being inspired with Hasimoto's 
pioneering work in \cite{Hasimoto}, 
Chang, Shatah and Uhlenbeck 
constructed a good moving frame 
along the map, 
and rigorously reduced the equation of 
the one-dimensional Schr\"odinger map 
into a compact Riemann surface 
to a simple form of a complex-valued 
nonlinear Schr\"odinger equation in \cite{CSU}. 
Using the same idea, the author studied 
the geometric reduction of the equations 
of higher-order dispersive flows in \cite{Onodera2}. 
In addition, 
time-global existence theorems were also studied 
under some geometric conditions. 
For the one-dimensional Schr\"odinger maps, 
time-global existence holds 
if $(N,J,g)$ is locally symmetric. 
See 
\cite{Koiso}, 
\cite{PWW2},
and 
\cite{SSB}. 
For the third-order equation \eqref{equation:pde}, 
Nishiyama and Tani in \cite{NT} and \cite{TN} proved 
time-local and time-global existence of solutions 
when $X=\mathbb{R}$ or $X=\mathbb{R}/\mathbb{Z}$, 
$N=\mathbb{S}^2$,  
and the integrability condition $b=a/2$ is satisfied. 
They made use of some conservation laws 
to prove the global existence theorem. 
These conservation laws were discovered by Zakharov and Shabat 
in the study of the Hirota equation. 
See \cite{ZS} for details. 
In \cite{Onodera1} the author generalized these results 
when $X=\mathbb{R}/\mathbb{Z}$. 
He proved a time-local existence theorem for
\eqref{equation:pde}-\eqref{equation:data}
when 
$N$ is a compact K\"ahler manifold, 
and proved a time-global existence theorem when 
$N$ is a compact Riemann surface 
with a constant curvature $K$, and 
the condition $b=Ka/2$ holds. 
\par
On the other hands, 
almost Hermitian manifolds 
do not necessarily satisfy 
the K\"ahler condition $\nabla{J}=0$. 
For example, it is well-known that 
$\mathbb{S}^6$, 
the Hopf manifold 
$\mathbb{S}^{2p+1}\times\mathbb{S}^{1}$, 
and 
$\mathbb{S}^{2p+1}\times\mathbb{S}^{2q+1}$ 
($p,q=1,2,3,\dotsc$) 
never admit the structure of K\"ahler manifolds. 
If the K\"ahler condition fails to hold, 
then $\nabla{J}$ causes the so-called loss of 
one-derivative, 
and the equation \eqref{equation:pde} behaves like 
the Cauchy-Riemann equation. 
In this case, 
the classical energy method breaks down. 
The main purpose of this paper is to show 
the time-local existence theorem of 
\eqref{equation:pde}-\eqref{equation:data} 
without the K\"ahler condition. 
To state our results, 
we here introduce some function spaces for mappings. 
\begin{definition} 
\label{definition:sobolev}
Let $\mathbb{N}$ be the set of positive integers. 
For $m\in\mathbb{N}\cup\{0\}$, 
the Sobolev space of mappings is defined by  
$$
 H^{m+1}(\RR;N)
 =
 \{
u \in C(\RR; N)
\ \vert \ 
u_x \in H^m(\RR;TN) 
\},
$$
where $u_x \in H^m(\RR;TN)$ means that $u_x$ satisfies 
$$
\lVert{u_x}\rVert_{H^{m}(\RR;TN)}^2
=
\sum_{j=0}^m
\int_{\RR} 
g_{u(x)}(\nabla_x^ju_x(x),\nabla_x^ju_x(x))dx
<
+\infty.
$$
Moreover, let $I$ be an interval in $\RR$, 
and let $w$ be an isometric embedding of 
$(N,J,g)$ into the standard Euclidean space 
$(\RR^d,g_0)$.
We say that $u\in C(I;H^{m+1}(\RR;N))$ 
if $u\in C(I\times \RR;N)$ and 
$(w{\circ}u)_x\in C(I;H^{m}(\RR;\RR^d))$,
where $C(I;H^{m}(\RR;\RR^d))$ is 
the set of usual Sobolev space valued continuous functions on $I$. 
\end{definition}
%
%
Our main results is the following.
%
\begin{theorem}
\label{theorem:eo}
Let $(N,J,g)$ be a compact almost Hermitian manifold, 
and let $a\neq 0$, $b\in \RR.$
Then for any $u_0{\in}H^{m+1}(\RR;N)$ 
with an integer $m\geqslant 4$, 
there exists a constant $T>0$ 
depending only on $a$, $b$, $N$ and 
$\lVert{u_{0x}}\rVert_{H^4(\RR;TN)}$ 
such that the initial value problem 
\eqref{equation:pde}-\eqref{equation:data} 
possesses a unique solution 
$u{\in}C([-T,T];H^{m+1}(\RR;N))$. 
\end{theorem}
Roughly speaking, Theorem~\ref{theorem:eo} says that 
\eqref{equation:pde}-\eqref{equation:data}
has a time-local solution in the usual Sobolev space 
$H^5(\RR;\RR^d)=(1-\p_x^2)^{-5/2}L^2(\RR;\RR^d)$. 
\par
Our idea of the proof comes from 
the theory of linear dispersive 
partial differential operators. 
Consider the initial value problem for 
linear partial differential equations of the form 
\begin{equation}
u_t+u_{xxx}+a(x)u_x+b(x)u=f(t,x)
\quad\text{in}\quad
\mathbb{R}\times\mathbb{R},
\label{equation:oogaki} 
\end{equation}
where 
$a(x),b(x)\in\mathscr{B}^{\infty}(\mathbb{R})$, 
which is the set of all smooth functions on $\mathbb{R}$ 
whose derivative of any order are bounded on $\mathbb{R}$, 
$u(t,x)$ is a complex-valued unknown function, 
and 
$f(t,x)$ is a given function. 
Tarama proved in \cite{Tarama} that 
the initial value problem for \eqref{equation:oogaki} 
is $L^2$-well-posed 
if and only if 
\begin{equation}
\left\lvert
\int_x^y\operatorname{Im}\hspace{1pt}a(s)ds
\right\rvert
\leqslant
C\lvert{x-y}\rvert^{1/2}
\label{equation:taramacondition} 
\end{equation}
for any $x,y\in\mathbb{R}$ 
with some constant $C>0$. 
The necessity is proved by 
the usual method of asymptotic solutions. 
In order to prove the sufficiency, 
Tarama first constructed 
a nice pseudodifferential operators of order zero 
which is automorphic on 
$L^2(\mathbb{R};\mathbb{C})$ 
under the condition 
\eqref{equation:taramacondition},
and eliminates 
$\sqrt{-1}\operatorname{Im}\hspace{1pt}a(x)\partial_x$. 
This is one of the methods of bringing out 
the local smoothing effect of $e^{-t\partial_x^3}$ on $\RR$, 
and this property breaks down on 
$\mathbb{R}/\mathbb{Z}$. 
See e.g., \cite{Doi}.  
Tarama also pointed out unofficially that if 
$\operatorname{Im}\hspace{1pt}a{\in}L^2(\mathbb{R};\mathbb{R})$, 
then 
\eqref{equation:taramacondition} 
holds 
and the proof of sufficiency becomes quite easier 
than the general case of 
\eqref{equation:taramacondition}.  
In this case, a gauge transformation defined by 
\begin{equation}
u(x)
\longmapsto
v(x)
=
u(x)
\exp\left(\frac{1}{3}\int_{-\infty}^x
\{\operatorname{Im}\hspace{1pt}a(y)\}^2
dy\right) 
\label{equation:kojimayoshio}
\end{equation}
is automorphic on $L^2(\mathbb{R};\mathbb{C})$, 
and \eqref{equation:oogaki} becomes 
\begin{equation}
v_t
+
v_{xxx}
-
\{\operatorname{Im}\hspace{1pt}a(x)\}^2
v_{xx}
+
\{\tilde{a}(x)
+
\sqrt{-1}\operatorname{Im}\hspace{1pt}a(x)\}
v_x
+
\tilde{b}(x)
v
=
\tilde{f}(t,x)
\label{equation:yuko}
\end{equation}
with some 
$\tilde{a},\tilde{b}\in\mathscr{B}^{\infty}(\mathbb{R})$ 
and $\tilde{f}$, where $\tilde{a}$ is a real-valued. 
The initial value problem for \eqref{equation:yuko} 
is $L^2$-well-posed in the positive direction of $t$ 
since the second-order term 
$\{\operatorname{Im}\hspace{1pt}a(x)\}^2
\partial_x^2$ 
dominates the seemingly bad first-order term 
$\sqrt{-1}\operatorname{Im}\hspace{1pt}a(x)
\partial_x$ 
essentially. 
In this special case, 
pseudodifferential calculus is not required. 
\par
We make use of the idea of the gauge transformation 
\eqref{equation:kojimayoshio}. 
Roughly speaking, we see $\nabla_x^mu_x$ satisfies the form 
\begin{equation}
\left(
\nabla_t-a\nabla_x^3-\nabla_xJ_u\nabla_x
\right)
\nabla_x^mu_x 
-m(\nabla_xJ_u)\nabla_x\nabla_x^mu_x
=
\text{harmless terms},
\label{equation:derivativeloss}
\end{equation}
where $(\nabla_xJ_u)$ is
the covariant derivative  
of the $(1,1)$-tensor  field $J_u$ with respect to $x$ along $u$.
The term $m(\nabla_xJ_u)\nabla_x\nabla_x^mu_x$ 
cannot be controlled by the classical energy method 
since 
$(\nabla_xJ_u)$
behaves as anti-symmetric operator on  
$L^2(\RR;TN)$ in the sense  
$$
\int_{\RR}g((\nabla_xJ_u)V, W)dx 
=
-\int_{\RR}g(V, (\nabla_xJ_u)W)dx, 
\quad
\text{for}
\quad
V,W \in \Gamma(u^{-1}TN).
$$
We introduce a gauge transformation on 
$u^{-1}TN$ defined by 
\begin{equation}
\nabla_x^mu_x(t,x)
\longmapsto
\nabla_x^mu_x(t,x)
\exp
\left(
-\frac{1}{3a}
\int_{-\infty}^x
g(u_x(t,y),u_x(t,y))
dy
\right),
\label{equation:mao} 
\end{equation}
which eliminates the bad term essentially 
since $(\nabla_xJ_u)=O\left( g(u_x,u_x)^{1/2}  \right)$. 
Parabolic regularization and 
the energy estimates with \eqref{equation:mao} 
prove Theorem~\ref{theorem:eo}. 
The assumption $m\geqslant 4$ is 
the requirement on the integer 
for our method to work.
\par
When $(N,J,g)$ is a K\"ahler manifold, 
we do not need the regularity $m\geqslant 4$. 
In this case, 
the term
$m(\nabla_xJ_u)\nabla_x\nabla_x^mu_x$ vanishes 
in \eqref{equation:derivativeloss}, thus the 
classical energy method works.
Indeed we prove the following.
\begin{theorem}
\label{theorem:eeo}
Let $(N,J,g)$ be a compact K\"ahler manifold and 
let $a\neq 0$ and $b\in \RR$. 
Then for any $u_0{\in}H^{m+1}(\RR;N)$ 
with an integer $m\geqslant2$, 
there exists a constant $T>0$ 
depending only on 
$a,b,N$, and $\lVert{u_{0x}}\rVert_{H^2(\RR;TN)}$ 
such that the initial value problem
\eqref{equation:pde}-\eqref{equation:data}
possesses a unique solution 
$u{\in}C([-T,T];H^{m+1}(\RR;N))$.
\end{theorem}
\begin{theorem}
\label{theorem:meo}
Let $(N,J,g)$ be a compact Riemann surface 
with constant Gaussian curvature $K$
and let $a\neq 0$ and $b=aK/2$. 
Then for any $u_0{\in}H^{m+1}(\RR;N)$ with an integer $m\geqslant2$, 
there exists a unique solution 
$u{\in}C(\RR;H^{m+1}(\RR;N))$
to
\eqref{equation:pde}-\eqref{equation:data}.
\end{theorem}
Theorem~\ref{theorem:eeo} and \ref{theorem:meo} 
are analogues of the results on $X=\RR/\mathbb{Z}$ 
in \cite{Onodera1}.
We remark that Theorem~\ref{theorem:meo} 
generalizes the results on $X=\RR$ in \cite{NT} and \cite{TN}. 
The key idea of the proof is the use of some conserved quantities 
generalizing what is used in  \cite{NT}.
Examples of Riemann surfaces satisfying the conditions in 
Theorem~\ref{theorem:meo} are not only 
the two-sphere $\mathbb{S}^2$ $(K=1)$ 
and  the flat torus 
$\mathbb{T}^2=\mathbb{R}^2/\mathbb{Z}^2$ 
($K=0$), 
but also closed hyperbolic surfaces ($K=-1$). 
\par 
The organization of this paper is as follows. 
Section \ref{section:note} is devoted to geometric preliminaries. 
In Section \ref{section:parabolic} 
we construct a sequence of approximate solutions 
by solving the IVP for a fourth-order parabolic equation. 
In Section \ref{section:energy} 
we obtain uniform estimates of approximate solutions. 
In Section \ref{section:proof1} 
we complete the proof of Theorem~\ref{theorem:eo}.
Finally, in Section \ref{section:special} 
we give the sketch of the proof of 
 Theorem~\ref{theorem:eeo} and \ref{theorem:meo}. 
\section{Geometric Preliminaries}
\label{section:note}
In this section, 
we introduce some geometric notations
used later in our proof. 
One can refer \cite{Nishikawa} 
for the elements of nonlinear 
geometric analysis.
\par
We will use 
$C=C(\cdot, \ldots, \cdot)$ 
to denote a positive constant depending on 
the certain parameters, 
geometric properties of $N$, 
et al.
The partial differentiation is written by $\p$, 
or the subscript, e.g., 
$\p_xf$, $f_x$, 
to distinguish from the covariant derivative 
along the curve, e.g., $\nabla_x$.
\par
Throughout this paper,   
$w$ is fixed as an isometric embedding mapping 
from $(N, J, g)$ into 
a standard Euclidean space $(\RR^d, g_0)$. 
Existence of $w$ is ensured by the celebrated works of 
Nash \cite{Nash}, 
Gromov and Rohlin \cite{GR}, 
and related papers. 
\par
For $\delta>0$, 
let $(w(N))_{\delta}$ be a $\delta$-tubular neighbourhood 
of  $w(N)\subset \RR^d$ defined by 
$$
(w(N))_{\delta}=
\left\{
Q=q+X \in \RR^d \ | \ 
q\in w(N), \
X\in (T_{q}w(N))^{\perp}, \ 
|X|< \delta \
\right\}
$$ 
where $|\cdot|$ denotes the distance in $\RR^d$,
and let $\pi :(w(N))_{\delta}\to w(N)$ 
be the nearest point projection map defined by 
$\pi(Q)=q$ for $Q=q+X\in (w(N))_{\delta}$. 
Since $w(N)$ is compact,  
for any sufficiently small $\delta$, 
$\pi$ exists and is smooth. 
We fix such small $\delta$. 
\par
Let $u:\RR\to N$ be given. 
$u^{-1}TN=\bigcup_{x\in \RR}T_{u(x)}N$ 
is the pull-back bundle induced 
from $TN$ by $u$. 
$V$ is called a section of $u^{-1}TN$ 
if $V(x)\in T_{u(x)}N$ for all $x\in \RR$.
We denote 
the space of all the sections of 
$u^{-1}TN$ 
by $\Gamma(u^{-1}TN)$.
For $V, W\in \Gamma(u^{-1}TN)$, 
define the quantities like $L^2$-inner product by
$$
\int_{\RR}
g(V, W)
dx
=
\int_{\RR}
g_{u(x)}(V(x),W(x))
dx, 
\quad 
\|V\|_{L^2(\RR;TN)}^2
=
\int_{\RR}
g(V, V)
dx.
$$
Then the quantity
$\|u_x\|_{H^m(\RR;TN)}^2$ 
defined in Definition~\ref{definition:sobolev} 
is written by
$$
\|u_x\|_{H^m(\RR;TN)}^2
=
\sum_{j=0}^m
\|\nabla_x^ju_x\|_{L^2(\RR;TN)}^2.
$$
In contrast, the standard $L^2$-product 
and $L^2$-norm  
are written by 
$$
\left\langle 
V,W
\right\rangle
=
\int_{\RR}
g_{0}(V(x),W(x))
dx, 
\quad
\|V\|_{L^2(\RR;\RR^d)}^2
=
\left\langle 
V,V
\right\rangle 
$$
for 
$V,W\in L^2(\RR;\RR^d),$
and the quantity
$\|V\|_{H^m(\RR;\RR^d)}^2$ 
is written by
$$
\|V\|_{H^m(\RR;\RR^d)}^2
=
\sum_{j=0}^m
\|\p_x^jV\|_{L^2(\RR;\RR^d)}^2.
$$
At this time 
$\|u_x\|_{H^m(\RR;TN)}<\infty$ 
if and only if 
$\|(w{\circ}u)_x\|_{H^m(\RR;\RR^d)}<\infty$. 
See, e.g., \cite[Section~1]{SZ} or \cite[Proposition~2.5]{MCGAHAGAN} 
for this equivalence.
Noting this equivalence, 
we see
$$
  H^{m+1}(\RR;N)
 =
 \{
u \in C(\RR; N)
\ \vert \ 
(w{\circ}u)_x\in H^{m}(\RR;\RR^d) 
\}.
$$
Finally, 
for $\alpha >0$, $m\in \mathbb{N}\cup \{0\}$ 
and an interval $I\subset \RR$,
$C^{0,\alpha}(I;H^m(\RR;\RR^d))$ 
denotes the usual $H^m(\RR;\RR^d)$-valued 
$\alpha$-H\"order space on $I$.
We will make use of fundamental Sobolev space theory 
of $H^{m}(\RR;\RR^d)$ later in our proof. 
\section{Parabolic Regularization}
\label{section:parabolic}
The aim of this section is 
to obtain a sequence 
$\{u^{\ep}\}_{\ep\in (0, 1)}$ 
solving 
\begin{alignat}{2}
  u_t
& = 
  -
  \ep\, \nabla_x^3u_x 
  +
  a\, \nabla_x^2u_x
  +
  J_u\nabla_xu_x
  +
  b\, g_u(u_x, u_x)u_x,
&
  \quad\text{in}\quad
& (0,T_{\ep})\times{\RR},
\label{equation:pde4}
\\
  u(0,x)
& =
  u_0(x)
&
  \quad\text{in}\quad
& {\RR}
\label{equation:data4}
\end{alignat}
for each $\ep\in (0, 1)$, 
where $u=u^{\ep}(t,x)$ 
is also an $N$-valued unknown function of 
$(t,x)\in [0, T_{\ep}]\times \RR$, 
and $u_0$ is the same initial data as that of 
\eqref{equation:pde}-\eqref{equation:data} 
independent of $\ep\in (0, 1)$.
The argument in this section is essentially 
same as that in  \cite[Section~3]{Onodera1}.
In fact, 
we can show that \eqref{equation:pde4}-\eqref{equation:data4} 
admits a unique solution near the initial data $u_0$. 
Define 
$$
L^{\infty}_{\delta,T }
=
\left\{
u\in L^{\infty}((0,T)\times \RR;N)
\ | \
\|
w{\circ}u-w{\circ}u_0
\|_{L^{\infty}((0,T)\times \RR;\RR^d)}
\leqslant \delta/2
\right\}
$$
for $T>0$, 
where $\delta >0$ is the fixed constant describing the radius of the tubular
neighbourhood of $w(N)$ as stated in the previous section. 
We show the following. 
%
%
\begin{proposition}
\label{proposition:pr}
Let $u_0\in H^{k+1}(\RR;N)$ with an integer $k\geqslant 2$. 
Then for each $\ep\in (0, 1)$, 
there exists a constant
$
T_{\ep}
=
T(\ep, a, b, N,
\| u_{0x}   \|_{H^k(\RR;TN)} )>0
$
and a unique solution 
$u=u^{\ep}\in C([0, T_{\ep}]; H^{k+1}(\RR; N))
\cap 
L^{\infty}_{\delta, T_{\ep}}$ 
to
\eqref{equation:pde4}-\eqref{equation:data4}. 
\end{proposition}
\begin{bew}{ of Proposition~\ref{proposition:pr}} 
%
Via the relation $v=w{\circ}u$, 
the IVP \eqref{equation:pde4}-\eqref{equation:data4} 
is equivalent to the following problem
\begin{alignat}{2}
  v_t
& = 
  -\ep v_{xxxx}
  +F(v)
&
  \quad\text{in}\quad
& (0,T_{\ep})\times{\RR},
\label{equation:pde5}
\\
  v(0,x)
& =
  w{\circ} u_0(x)
&
  \quad\text{in}\quad
& {\RR},
\label{equation:data5}
\end{alignat}
where 
$v=v^{\ep}(t,x)$ 
is a $w(N)$-valued unknown function of 
$(t,x)\in [0, T_{\ep}]\times \RR$, 
and $F(v)$ is written by the form
\begin{align}
 F(v) 
& =
 -\ep
  \{ [A(v)(v_x, v_x)]_{xx}
     +[A(v)(v_{xx}+A(v)(v_x, v_x), v_x)]_x
\nonumber
\\
& \qquad 
    +A(v)(v_{xxx}+[A(v)(v_x,v_x)]_x
    + A(v)(v_{xx}+A(v)(v_x, v_x), v_x), v_x)
   \}
\nonumber
\\
& \quad + a
  \{ v_{xxx} + [A(v)(v_x, v_x)]_x 
    + A(v)(v_{xx}+A(v)(v_x, v_x), v_x) 
  \}
\nonumber
\\
& \quad 
 + dw_{w^{-1}\circ v}J_{w^{-1}\circ v}dw^{-1}_v(v_{xx}+A(v)(v_x, v_x))
 + b|v_x|^2 v_x,
\nonumber 
\end{align}
where,
$
A(v)(\cdot, \cdot ): 
T_vw(N)\times T_vw(N)\to (T_vw(N))^{\perp}
$
is the second fundamental form 
of $w(N)\subset \RR^d$ 
at $v\in w(N)$.
Note that there exists 
$G\in C^{\infty}(\RR^{4d};\RR^d)$ 
such that
\begin{equation}
F(v)=G(v, v_x, v_{xx}, v_{xxx})
\nonumber
\end{equation}
for $v:\RR\to w(N)$, 
and
$G(v,p,q,r)$ satisfies 
$$
G(v, 0, 0, 0)=0,
\quad
\dfrac{\p^2 G}{\p r^2}(v,p,q,r)=0.
$$
The equation
\eqref{equation:pde5} 
is a system of fourth-order parabolic evolution equations 
for $\RR^d$-valued function. 
In place of the IVP \eqref{equation:pde4}-\eqref{equation:data4}, 
we will solve the IVP \eqref{equation:pde5}-\eqref{equation:data5}.
%
The proof consists of 
the following two steps.
First, we construct a solution of 
\eqref{equation:pde5}-\eqref{equation:data5} 
whose image are contained in 
$(w(N))_{\delta}\subset \RR^d$.
More precisely, we extend \eqref{equation:pde5} 
to an equation
for the vector-valued function valued in  $(w(N))_{\delta}$ 
and 
construct a unique time-local solution of the 
IVP for the 
extended equation 
in the class
$$
Y_T
 =
  \{
   v\in 
   X_T
   \ | \
   \| v-w{\circ} u_0 \|_{L^{\infty}((0,T)\times \RR;\RR^d)}
   \leqslant 
   \delta/2
  \}
$$
for sufficiently small $T>0$. Here 
$$
X_T 
=
\{
   v\in 
   C([0, T]\times \RR; \RR^d)
   \ | \  
   v_x\in C([0, T];H^{k}(\RR;\RR^d))
  \}
$$
is the Banach space with the following norm 
$$
\left\| v  \right\|_{X_T}
=
\left\| v  \right\|_{L^{\infty}([0,T]\times \RR;\RR^d)}
+
\left\| v_x  \right\|_{L^{\infty}(0,T; H^k(\RR; \RR^d))},
\quad
v\in X_T.
$$
Secondly, we check that this solution is actually $w(N)$-valued 
by using a kind of maximum principle. 
\par 
In short, it suffices to show the following two lemmas to complete 
our proof.
%
%
\begin{lemma}
\label{lemma:calculus}
For each $\varepsilon\in (0,1)$, 
there exists a constant
$
T_{\ep}>0
$
depending on 
$
\ep, a, b, N 
$
and
$
\| (w{\circ} u_0)_x   \|_{H^k(\RR;\RR^d)} 
$
and there exists a unique solution 
$v=v^{\ep}\in Y_{T_{\ep}}$
to
\begin{alignat}{2}
  v_t
& = 
  -\ep v_{xxxx}+F(\pi{\circ} v)
&
  \quad\mathrm{in}\quad
& (0,T_{\ep})\times{\RR},
\label{equation:pde6}
\\
  v(0,x)
& =
  w{\circ} u_0(x)
&
  \quad\mathrm{in}\quad
& {\RR}.
\label{equation:data6}
\end{alignat}
Moreover, the map 
$(w{\circ}u_0)_{x}\in H^{k}(\RR;\RR^d)\to 
v^{\ep}_x\in C([0,T_{\ep}];H^{k}(\RR;\RR^d))$ 
is continuous.
\end{lemma}
%
%
\begin{lemma}
\label{lemma:maximum}
Fix $\varepsilon\in (0,1)$. 
Assume that 
$v=v^{\ep}\in Y_{T_{\ep}}$ 
solves \eqref{equation:pde6}-\eqref{equation:data6}. 
Then $v(t,x)\in w(N)$ 
for all 
$(t,x)\in [0,T_{\ep}]\times \RR$, 
thus $v$ solves \eqref{equation:pde5}-\eqref{equation:data5}. 
\end{lemma}
\begin{bew}{of Lemma~\ref{lemma:calculus}}
The idea of the proof is due to the contraction mapping argument. 
\par
Let $L$ be a nonlinear map defined by 
\begin{align}
Lv(t)
&=
e^{-\ep t\p_x^4}v_0 
+
\int_0^te^{-\ep (t-s)\p_x^4}F((\pi{\circ} v)(s))ds
\nonumber 
\\
&=
\int_{\RR}
E(t,x-y)v_0(y)dy
+
\int_0^t
\int_{\RR}
E(t-s, x-y) 
F((\pi{\circ} v)(s,y))
dyds, 
\nonumber
\end{align}
where  $v_0=w{\circ}u_0$,
and 
$E(t,x)$
is the fundamental solution associated to 
$\p_t+\ep \p_x^4$. 
Note that,
if $v\in Y_T$, 
$\pi {\circ} v$ takes value in $w(N)$ and thus 
$F(\pi {\circ} v)$ makes sense.
The IVP 
\eqref{equation:pde6}-\eqref{equation:data6} 
is equivalent to an integral equation of the form 
$v=Lv$.
\par
Set $M=\|v_{0x}\|_{H^k(\RR;\RR^d)}$, 
and 
define the space
$$
Z_T
 =
  \{
   v\in 
   Y_T \ | \
   \| v_x  \|_{L^{\infty}(0, T; H^k(\RR;\RR^d))}
   \leqslant
   2M
  \}.
$$
$Z_T$ is a closed subset of the Banach space $X_T$. 
To complete the proof,
we have only to show that the map $L$ has a unique fixed point 
in $Z_{T_{\ep}}$ for sufficiently small $T_{\ep}>0$, 
since the uniqueness in the whole space 
$Y_{T_{\ep}}$ follows by similar and standard arguments. 
\par
First, consider the properties of $e^{-\ep t\p_x^4}$. 
Since $u_0\in H^{k+1}(\RR;N)$, 
$v_0$ is especially bounded 
and uniformly continuous on $\RR$. 
Thus, it is easy to check that
\begin{equation}
e^{-\ep t\p_x^4}v_0
\longrightarrow 
v_0
\quad 
\text{in} 
\quad 
C(\RR;\RR^d) 
\quad 
\text{as} 
\quad 
t\to 0, 
\label{equation:fundamental0}
\end{equation}
and
\begin{equation} 
\| e^{-\ep t\p_x^4}v_{0x} \|_{H^{k}(\RR;\RR^d)}
 \leqslant 
  \|v_{0x} \|_{H^{k}(\RR;\RR^d)}.
\label{equation:fundamental}
\end{equation}
Moreover,
 $e^{-\ep t\p_x^4}$ 
gains the regularity of order $3$, 
since 
$
(\ep^{1/4}t^{1/4}|\xi|)^j
e^{-\ep t\xi^4}
$
is bounded 
for $j=0, 1, 2, 3$. 
In fact, there exists $C_1>0$ 
such that
\begin{equation}
\| e^{-\ep t\p_x^4}\phi \|_{H^{k+1}(\RR;\RR^d)}
 \leqslant
  C_1 \ep^{-3/4}t^{-3/4}
  \|\phi \|_{H^{k-2}(\RR;\RR^d)}
\label{equation:smoothing}
\end{equation}
holds for any 
$\phi \in H^{k-2}(\RR;\RR^d)$.
\par
Secondly, consider the nonlinear estimates of 
$F(\pi{\circ} v)$.
If $v$ belongs to the class $Z_{T}$, 
we see $v(t, \cdot)\in C(\RR; (w(N))_{\delta})$ 
and 
$\|v_x(t)\|_{H^k(\RR;\RR^d)}
\leqslant 
2M$
follows for all $t\in [0, T]$. 
Thus, by observing the form of $F(v)$ and the compactness of $w(N)$, 
it is easy to check that there exists 
$C_2=C_2(a,b,M,N)>0$ such that 
\begin{align}
  \| F(\pi{\circ} v)(t)\|_{H^{k-2}(\RR;\RR^d)}
 &\leqslant
  C_2\| v_x(t) \|_{H^k(\RR;\RR^d)},
\label{equation:composite1}
\\
  \| F(\pi{\circ} u)(t)-F(\pi{\circ} v)(t)\|_{H^{k-2}(\RR;\RR^d)}
 &\leqslant
  C_2
 \left(  
 \| u(t)-v(t) \|_{L^{\infty}(\RR;\RR^d)} 
  +    
 \| u_x(t)-v_x(t) \|_{H^k(\RR;\RR^d)}
 \right)
\label{equation:composite2}
\end{align}
for any 
$u, v\in Z_T$. 
\par 
Using the properties 
\eqref{equation:fundamental0},
\eqref{equation:fundamental},  
\eqref{equation:smoothing} 
and the nonlinear estimates 
\eqref{equation:composite1}, 
\eqref{equation:composite2}, 
we can prove that $L$ is a contraction mapping 
from $Z_{T_{\ep}}$ into itself 
if $T_{\ep}$ is sufficiently small. 
It is the standard argument, 
thus we omit the rest of the proof.
\qed
\end{bew}
\begin{remark}
\label{remark:gainofparabolic}
Suppose that $v^{\ep}\in Y_{T_{\ep}}$ solves  
\eqref{equation:pde5}-\eqref{equation:data5}. 
Then we can easily check 
$v^{\ep}_{xxxx}\in L^2(0, T_{\ep}; H^1(\RR;\RR^d))$ 
and 
$F(\pi\circ v^{\ep})\in L^2(0, T_{\ep}; H^1(\RR;\RR^d))$ 
from the standard arguments. 
Thus we see
$v^{\ep}_t$ belongs to the same class $L^2(0, T_{\ep}; H^1(\RR;\RR^d))$, 
which implies that 
$v^{\ep}-v_0$ belongs to the class 
$C^{0,1/2}([0,T_{\ep}];H^1(\RR;\RR^d))$. 
\end{remark}
%
\begin{bew}{of Lemma~\ref{lemma:maximum}}
Suppose 
$v\in Y_{T_{\ep}}$ 
solves \eqref{equation:pde6}-\eqref{equation:data6}.
Define the map $\rho: (w(N))_{\delta}\to \RR^d$ by 
$\rho(Q)=Q-\pi(Q)$ for $Q\in (w(N))_{\delta}$. 
Then we deduce
\begin{align}
\left|
\rho {\circ} v(t,x)
\right|
&=
\min_{q\in w(N)}
\left|
v(t,x)-q
\right|
\leqslant 
|v(t,x)-v_0(x)|.
\nonumber 
\end{align}
Notice that the first equality above is due to 
the compactness of $w(N)$.  
In addition, 
as is stated in Remark~\ref{remark:gainofparabolic},
$v(t)-v_0$ 
belongs to 
$L^{2}(\RR;\RR^d)$ and 
thus $\rho {\circ} v(t)$ 
makes sense in $L^2(\RR;\RR^d)$ for each $t$. 
To obtain that $v$ is $w(N)$-valued, we will show 
$$
\| \rho {\circ} v(t)\|_{L^2(\RR;\RR^d)}^2
=
\left\langle
 \rho {\circ} v(t), 
  \rho {\circ} v(t) 
\right\rangle
=0
$$
for all $t\in[0,T_{\ep}]$.
Since $\pi + \rho$ is identity on $(w(N))_{\delta}$,
\begin{equation}
 d\pi_v + d\rho_v
 =
 I_d 
\label{equation:decomposition}
\end{equation}
holds on $T_v(w(N))_{\delta}$, 
where $I_d$ is the identity. 
By identifying 
$T_v(w(N))_{\delta}$ 
with 
$\RR^d$, 
we see that
$v_t(t,x)\in T_{v(t,x)}(w(N))_{\delta}$
and
$
d\pi_{v}(v_t)(t,x)
\in 
T_{\pi\circ v(t,x)}w(N)
$ 
for each $(t,x)$. 
Thus it follows that 
$
\left\langle
\rho {\circ} v, d\pi_v ( v_t) 
\right\rangle
=0.
$
Using this relation and \eqref{equation:decomposition}, we deduce
$$
\frac{1}{2}\frac{d}{dt}
 \| \rho {\circ} v  \|_{L^2(\RR;\RR^d)}^2
=
\left\langle
  \rho {\circ} v, d\rho_v ( v_t) 
 \right\rangle
=
 \left\langle
  \rho {\circ} v, d\rho_v ( v_t)+d\pi_v ( v_t)
 \right\rangle
=
 \left\langle
  \rho {\circ} v, v_t 
 \right\rangle.
$$
Recall here, 
by the form of the right hand side of 
\eqref{equation:pde5},
that 
$-\ep \tilde{v}_{xxxx}+F(\tilde{v})
\in 
\Gamma(\tilde{v}^{-1}Tw(N))$ 
holds for any 
$\tilde{v}:\RR\to w(N)$.
Thus we see 
$(-\ep (\pi{\circ} v)_{xxxx}+F(\pi{\circ} v))(t)
\in 
\Gamma((\pi{\circ} v(t))^{-1}Tw(N))$ 
since $\pi{\circ} v(t)\in w(N)$, 
and thus this is perpendicular to 
$\rho{\circ} v(t)$. 
Noting this and substituting \eqref{equation:pde6}, we get
\begin{alignat}{2}
 \frac{1}{2}
 \frac{d}{dt}
 \| \rho {\circ} v  \|_{L^2(\RR;\RR^d)}^2
& = 
 \left\langle
 \rho {\circ} v, 
  -\ep v_{xxxx}
  +F(\pi{\circ} v) 
 \right\rangle
\nonumber
\\ 
&=
 \left\langle
 \rho {\circ} v, 
  - \ep(\rho{\circ} v)_{xxxx}
  - \ep(\pi{\circ} v)_{xxxx}
  +F(\pi{\circ} v) 
 \right\rangle
\nonumber
\\
& =
 \left\langle
 \rho {\circ} v, 
 - \ep(\rho{\circ} v)_{xxxx} 
 \right\rangle
\nonumber \\
&=
  - \ep
 \| (\rho {\circ} v)_{xx} \|_{L^2(\RR;\RR^d)}^2
\leqslant 
0,
\nonumber
\end{alignat}
which implies 
$
 \| \rho {\circ} v(t)\|_{L^2(\RR;\RR^d)}^2
 \leqslant
 \| \rho {\circ} v_0  \|_{L^2(\RR;\RR^d)}^2
 =0.
$ 
Hence 
$\rho {\circ} v(t)\equiv 0$ 
holds. 
Thus $v(t)$ is $w(N)$-valued for all $t$, 
which completes the proof.
\qed
\end{bew}
Set $u=w^{-1}{\circ}v$ for the solution $v$ 
in Lemma~\ref{lemma:calculus}. 
It is now obvious that this $u$ solves 
\eqref{equation:pde4}-\eqref{equation:data4}. 
Thus we complete the proof.
\qed
\end{bew}
\section{Geometric  energy estimates}
\label{section:energy}
Let $\{u^{\ep}\}_{\ep\in (0, 1)}$ 
be a sequence of solutions to 
\eqref{equation:pde4}-\eqref{equation:data4} 
constructed in Section~\ref{section:parabolic} 
with $k=m\geqslant 4$.  
We will  obtain the uniform estimate of 
$\{u_x^{\ep}\}_{\ep\in (0, 1)}$ and the existence time.
Our goal of this section is the following.
\begin{lemma}
\label{lemma:Energy} 
Let 
$u_0\in H^{m+1}(\RR;N)$ 
with an integer 
$m\geqslant 4$, 
and let $\{u^{\ep}\}_{\ep\in (0, 1)}$ 
be a sequence of solutions to 
\eqref{equation:pde4}-\eqref{equation:data4}.  
Then there exists a constant  
$T>0$ 
depending only on 
$a, b, N, \|u_{0x}\|_{H^4(\RR;TN)}$  
such that 
$\{u_x^{\ep}\}_{\ep\in (0, 1)}$ 
is a bounded sequence in 
$L^{\infty}(0, T; H^m(\RR;TN))$.
\end{lemma}
\begin{bew}{of Lemma~\ref{lemma:Energy}}
We define 
\begin{align}
K^{\ep}(t,x)
&=
-\frac{1}{3a}
\int_{-\infty}^x
g\left( 
u^{\ep}_x(t,y),u^{\ep}_x(t,y)
\right)
dy,
\nonumber 
\\
V^{\ep, (m)}(t,x)
&=
e^{K^{\ep}(t,x)}\nabla_x^mu_x^{\ep}(t,x),
\nonumber 
\\
N^{\ep}_m(t)
&=
\left(
 \| u_{x}^{\ep}(t) \|_{H^{m-1}(\RR;TN)}^2
+
 \| V^{\ep, (m)}(t) \|_{L^2(\RR;TN)}^2
\right)^{1/2}.
\nonumber
\end{align}
We will obtain the differential inequality for 
$\left(  N^{\ep}_m(t)  \right)^2$. 
Since $N^{\ep}_4(0)$ is independent of $\ep$, 
we set 
$r_0=N^{\ep}_4(0)$ and 
$$
T_{\ep}^{*}=
\sup
\left\{
T>0 \ | \ 
N^{\ep}_4(t)\leqslant 2r_0 
\ 
\text{for all}
\
t\in [0,T]
\right\}.
$$
Lemma~\ref{lemma:calculus} shows 
$T_{\ep}^{*}>0$.
Moreover, there exists a positive constant 
$C(a,r_0)>1$ such that 
$$
C(a,r_0)^{-1}N^{\ep}_m(t)
\leqslant 
\| u_{x}^{\ep}(t) \|_{H^{m}(\RR;TN)}
\leqslant 
C(a,r_0)N^{\ep}_m(t)
\quad 
\text{for $t\in [0, T_{\ep}^{*}]$.}
$$ 
This follows from the relation 
$$
\left| e^{\pm K^{\ep}(t,x)} \right|
\leqslant 
1+e^{ \frac{1}{3|a|} \| u_{x}^{\ep}(t) \|_{L^2(\RR;TN)}^2}
\leqslant 
1+e^{ \frac{1}{3|a|} \| u_{0x} \|_{L^2(\RR;TN)}^2}.
$$
Note here that the second inequality of the 
 estimate above is due to 
\begin{equation}
\| u_{x}^{\ep}(t) \|_{L^2(\RR;TN)}^2
\leqslant 
\| u_{0x} \|_{L^2(\RR;TN)}^2,
\nonumber
\end{equation}
which follows from the energy inequality of the form 
\begin{align}
\frac{1}{2}
\|u_{x}^{\ep} \|_{L^2(\RR;TN)}^2 
&=
\int_{\RR}
g\left( 
\nabla_tu_x^{\ep},u_x^{\ep}
\right)dx
\nonumber 
\\
&=
\int_{\RR}
g\left( 
\nabla_xu_t^{\ep},u_x^{\ep}
\right)dx
\nonumber
\\
&=\int_{\RR}
g\left(
-\ep\nabla_x^4u_x^{\ep}
+
a\nabla_x^3u_x^{\ep}
+
\nabla_xJ_{u^{\ep}}\nabla_xu_x^{\ep}
+
b\nabla_x[g(u_x^{\ep},u_x^{\ep})u_x^{\ep}]
, 
u_x^{\ep}
\right)dx
\nonumber 
\\
&=
-\ep\|\nabla_x^2u_x^{\ep} \|_{L^2(\RR;TN)}^2\leqslant 0.
\nonumber 
\end{align}
The last equality of the estimate above 
is easily checked by repeatedly using integration by parts. 
Especially, we see that  
$$
\int_{\RR}
g\left(
\nabla_xJ_{u^{\ep}}\nabla_xu_x^{\ep}
, 
u_x^{\ep}
\right)dx
=
-\int_{\RR}
g\left(
J_{u^{\ep}}\nabla_xu_x^{\ep}
, 
\nabla_xu_x^{\ep}
\right)dx
=0,
$$
where the second equality above is due to 
the fact that 
$(N,J,g)$ is an almost hermitian manifold.
\par
Having these notations and properties in mind, 
we show the following. 
\begin{proposition}
\label{proposition:Energyestimate}
There exists a positive constant 
$C=C(a,b,m,N,r_0)>0$ 
and an increasing function $P(\cdot)$ on $[0,+\infty)$ 
such that 
\begin{equation}
\label{equation:Energyinequality}
\begin{aligned}
&\frac{1}{2}\frac{d}{dt}
\left(
N^{\ep}_m(t)
\right)^2
+
\frac{\ep}{2}
\left(
\|\nabla_x^2 V^{\ep, (m)}(t) \|_{L^2(\RR;TN)}^2
+
\sum_{l=0}^{m-1}
\|\nabla_x^{l+2}u_x^{\ep}(t) \|_{L^2(\RR;TN)}^2
\right)
\\
&\qquad +
\frac{1}{2}
\|\left( g(u_x^{\ep}(t),u_x^{\ep}(t))\right)^{1/2}
\nabla_xV^{\ep, (m)} (t) \|_{L^2(\RR;TN)}^2 
\\
&\leqslant 
C(a,b,m,N,r_0)
P(N^{\ep}_4(t)+N^{\ep}_{m-1}(t))
\left(
N^{\ep}_m(t)
\right)^2
\end{aligned}
\end{equation}
follows for all $t\in [0,T_{\ep}^{*}]$.
\end{proposition}
\begin{bew}{of Proposition~\ref{proposition:Energyestimate}}
Throughout the proof of \eqref{equation:Energyinequality} 
we simply write 
$u$, $J$, $g$, $K$, $V^{(m)}$ 
in place of 
$u^{\ep}$, $J_{u^{\ep}}$, $g_{u^{\ep}}$, $K^{\ep}$, $V^{\ep,(m)}$  
respectively, 
and write 
$\|\cdot\|_{H^k}=\|\cdot\|_{H^k(\RR;TN)}$, 
$\|\cdot\|_{L^2}=\|\cdot\|_{L^2(\RR;TN)}$, 
$\|\cdot\|_{L^{\infty}}=\|\cdot\|_{L^{\infty}(\RR;TN)}$ 
for $k\in \mathbb{N}$,
and sometimes omit to write time variable $t$. 
\par
The main object of the proof is the estimation of 
\begin{equation}
\frac{1}{2}\frac{d}{dt}
\|V^{(m)}(t)\|_{L^2}^2
=
\int_{\RR}
g( \nabla_tV^{(m)}(t),  V^{(m)}(t))dx.
\label{equation:E1}
\end{equation}
Thus let us compute the equation of $V^{(m)}$. 
Operating $e^K\nabla_x^{m+1}$ on \eqref{equation:pde4}, 
we have
\begin{equation}
\nabla_tV^{(m)}
+
\ep \nabla_x^4V^{(m)}
-a\nabla_x^3V^{(m)}
-\nabla_xJ\nabla_xV^{(m)}
-\ep F_1
-F_2
=
F_3, 
\label{equation:Equation}
\end{equation}
where
\begin{align}
F_1
=&
4K_x\nabla_x^3V^{(m)}
+
6(K_{xx}-K_x^2)\nabla_x^2V^{(m)} 
\nonumber 
\\
&\quad +
4(K_{xxx}-3K_xK_{xx}+K_x^3)\nabla_xV^{(m)}
\nonumber 
\\
&\quad \quad +
(K_{xxxx}-4K_xK_{xxx}-3K_{xx}^2+6K_x^2K_{xx}-K_x^4)V^{(m)}
\nonumber 
\\
&\quad \quad \quad +
\sum_{l=0}^{m-1}e^K\nabla_x^l
\left[
R(u_x, \nabla_x^3u_x)
\nabla_x^{m-1-l}u_x
\right],
\label{equation:F1}
\\
F_2
=&
-3aK_x\nabla_x^2V^{(m)}
-3a(K_{xx}-K_x^2)\nabla_xV^{(m)}
\nonumber 
\\
&\quad -K_x\nabla_xJV^{(m)}-K_xJ\nabla_xV^{(m)}
+m\left( \nabla_xJ  \right)\nabla_xV^{(m)}
\nonumber 
\\
&\quad \quad -aR(u_x, \nabla_xV^{(m)} )u_x
+2b\, g(\nabla_xV^{(m)},u_x)u_x 
+b\, g(u_x,u_x)\nabla_xV^{(m)},
\label{equation:F2}
\\
F_3
=&
K_tV^{(m)}
-a\left(
\sum_{l=0}^{m-1}
e^K\nabla_x^l
\left[
R(u_x, \nabla_x^2u_x)\nabla_x^{m-1-l}u_x
\right]
-R(u_x, \nabla_xV^{(m)})u_x
\right)
\nonumber 
\\
&\quad
-\sum_{l=0}^{m-1}
e^K\nabla_x^l
\left[
R(u_x, J\nabla_xu_x)\nabla_x^{m-1-l}u_x
\right] 
\nonumber 
\\
&\quad \quad
-a(K_{xxx}-3K_{x}K_{xx}+K_x^3)V^{(m)}
\nonumber 
\\
&\quad \quad \quad
-(K_{xx}-K_x^2)JV^{(m)}
-mK_x\left(\nabla_xJ \right)V^{(m)} 
\nonumber 
\\
&\quad \quad \quad \quad
+\sum_{l=1}^m\sum_{j=1}^l
 \frac{l!}{j!(l-j)!}
 e^K\left(\nabla_x^{j+1}J\right)\nabla_x^{m+1-j}u_x
\nonumber 
\\
&\quad \quad \quad \quad \quad
-2bK_xg(V^{(m)},u_x)u_x-bK_xg(u_x,u_x)V^{(m)}
\nonumber 
\\
&\quad \quad \quad \quad \quad \quad
+
b\, e^K
       \sum_{\begin{smallmatrix}
        \alpha+\beta+\gamma=m+1 \\
        \alpha, \beta, \gamma\geqslant 0\\
        \max\{\alpha, \beta, \gamma\}\leqslant m
       \end{smallmatrix}}
  \frac{(m+1)!}{\alpha!\beta!\gamma!}
  g(
  \nabla_x^{\alpha}u_x, \nabla_x^{\beta}u_x
  )
  \nabla_x^{\gamma}u_x.
\label{equation:F3} 
\end{align}
Here $R$ denotes the curvature tensor on $(N, J, g)$, 
and 
$\left(\nabla_xJ\right)$ is the covariant derivative of 
$(1,1)$-tensor field $J$ with respect to $x$ along $u$ 
defined as 
\begin{equation}
\left(\nabla_xJ\right)V
=
\nabla_xJV-J\nabla_xV
\quad
\text{for}
\quad
V\in \Gamma(u^{-1}TN).
\label{equation:tensorJ0}
\end{equation} 
$\left(\nabla_xJ\right)$
is, by definition, a $(1,1)$-tensor field.
In the same way, $(\nabla_x^{j+1}J )$
denoting the $(j+1)$-th covariant derivative of $J$
is also $(1,1)$-tensor field along $u$.
See, Appendix, for the precise computations above.
\par
We next obtain the estimate of \eqref{equation:E1}
by putting \eqref{equation:Equation} into there. 
To make this estimate be clear 
or to focus only on the estimation of important 
parts as possible,
we use the notation as follows. 
\begin{definition}
For $A,B\in \RR$, 
$A\equiv B$ if and only if there exists a positive constant 
$C=C(a,b,m,N,r_0)>0$ 
and an increasing function $P(\cdot)$ on $[0,+\infty)$ 
such that
$$
A-B
\leqslant 
C(a,b,m,N,r_0)
P(N^{\ep}_4(t)+N^{\ep}_{m-1}(t))
\left(
N^{\ep}_{m}(t)
\right)^2
$$
follows for $t\in [0,T_{\ep}^{*}]$.
\end{definition}
First, it follows from the repeatedly using of
 integration by parts that
\begin{align}
\int_{\RR}
g(
-\ep\nabla_x^4V^{(m)},
V^{(m)}
)
dx
&=
-\ep
\|
\nabla_x^2V^{(m)}
\|_{L^2}^2,
\label{equation:EE1}
\\
\int_{\RR}
g(
a\nabla_x^3V^{(m)},
V^{(m)}
)
dx
&=
-a\int_{\RR}
g(
\nabla_x^2V^{(m)},
\nabla_xV^{(m)}
)
dx
=0,
\label{equation:EE2}
\\
\int_{\RR}
g(
\nabla_xJ\nabla_xV^{(m)},
V^{(m)}
)
dx
&=
-\int_{\RR}
g(
J\nabla_xV^{(m)},
\nabla_xV^{(m)}
)
dx
=0.
\label{equation:EE3}
\end{align}
\par
Next, let us go to the estimation of $F_2$. 
The following four terms 
\begin{gather*}
-3a(K_{xx}-K_x^2)\nabla_xV^{(m)}, 
\quad 
-aR(u_x,\nabla_xV^{(m)})u_x,
\\
2b\,g(\nabla_xV^{(m)}, u_x)u_x, 
\quad
b\,g(u_x, u_x)\nabla_xV^{(m)}
\end{gather*}
are easily controlled by a use of integration by parts. 
Indeed, we have 
\begin{equation}
\begin{aligned}
&\int_{\RR}
g(
-3a(K_{xx}-K_x^2)\nabla_xV^{(m)},
V^{(m)}
)dx  
\\
&=
-\frac{3a}{2}
\int_{\RR}
g(
(K_{xx}-K_x^2)\nabla_xV^{(m)},
V^{(m)}
)dx 
\\
&\quad \quad +
\frac{3a}{2}
\int_{\RR}
g(
(K_{xx}-K_x^2)V^{(m)},
\nabla_xV^{(m)}
)dx 
\\
&\quad \quad \quad +
\frac{3a}{2}
\int_{\RR}
g(
(K_{xx}-K_x^2)_xV^{(m)},
V^{(m)}
)dx 
\\
&=
\frac{3a}{2}
\int_{\RR}
g(
(K_{xx}-K_x^2)_xV^{(m)},
V^{(m)}
)dx 
\\
&\equiv  0,
\end{aligned}
\label{equation:E3}
\end{equation}
\begin{equation}
\begin{aligned}
&\int_{\RR}
g(
-aR(u_x, \nabla_xV^{(m)})u_x,
V^{(m)}
)dx 
\\
&=
-\frac{a}{2}
\int_{\RR}
g(
R(u_x, \nabla_xV^{(m)})u_x,
V^{(m)}
)dx  
\\
&\quad \quad +
\frac{a}{2}
\int_{\RR}
g(
R(u_x, V^{(m)})u_x,
\nabla_xV^{(m)}
)dx 
\\
&\quad \quad \quad +
\frac{a}{2}
\int_{\RR}
g(
R(u_x, V^{(m)})\nabla_xu_x,
V^{(m)}
)dx  
\\
&\quad \quad \quad \quad +
\frac{a}{2}
\int_{\RR}
g(
R(\nabla_xu_x, V^{(m)})u_x,
V^{(m)}
)dx 
\\
&\quad \quad \quad \quad \quad +
\frac{a}{2}
\int_{\RR}
g(
(\nabla_xR)(u_x, V^{(m)})u_x,
V^{(m)}
)dx
\\
&=
a
\int_{\RR}
g(
R(u_x, V^{(m)})\nabla_xu_x,
V^{(m)}
)dx  
\\
&\quad \quad +
\frac{a}{2}
\int_{\RR}
g(
(\nabla_xR)(u_x, V^{(m)})u_x,
V^{(m)}
)dx 
\\
&\equiv 0,
\end{aligned}
\label{equation:E4}
\end{equation}
\begin{equation}
\begin{aligned}
&\int_{\RR}
g(
2b\,g(\nabla_xV^{(m)},u_x )u_x,
V^{(m)}
)dx  
\\
&=-2b
\int_{\RR}
g(
g (V^{(m)},\nabla_xu_x )u_x,
V^{(m)}
)dx 
\\
&\equiv 0
\end{aligned}
\label{equation:E5}
\end{equation}
\begin{equation}
\begin{aligned}
&\int_{\RR}
g(
b\, g (u_x,u_x )\nabla_xV^{(m)},
V^{(m)}
)dx  
\\
&=-b
\int_{\RR}
g(
g( \nabla_xu_x,u_x )V^{(m)},
V^{(m)}
)dx 
\\
&\equiv 0.
\end{aligned}
\label{equation:E6}
\end{equation}
Notice that the second equality of \eqref{equation:E4} 
follows from the fundamental property of 
the Riemannian curvature tensor $R$ such as 
$$g\left( R(X,Y)Z,W  \right)=g\left( R(Z,W)X,Y  \right)
\quad 
\text{for}
\quad
X,Y,Z,W\in \Gamma(u^{-1}TN).
$$
\par
The estimates of the  rest terms of $F_2$ are demonstrated as follows.
For the estimate related to the term $-3aK_x\nabla_x^2V^{(m)}$, we have 
\begin{equation}
\begin{aligned}
&\int_{\RR}
g(
-3aK_x\nabla_x^2V^{(m)},
V^{(m)}
)dx 
\\
&=
\int_{\RR}
g(
g(u_x,u_x)\nabla_x^2V^{(m)},
V^{(m)}
)dx 
\\
&=
-\int_{\RR}
g(
g(u_x,u_x)\nabla_xV^{(m)}, \nabla_xV^{(m)}
)dx 
 -
2\int_{\RR}
g(
g(\nabla_xu_x,u_x)\nabla_xV^{(m)}, V^{(m)}
)dx 
\\
&=
-\|
\left(
g(u_x,u_x)
\right)^{1/2}
\nabla_xV^{(m)}
\|_{L^2}^2
+
\int_{\RR}
g(
\left[
g(\nabla_xu_x,u_x) 
\right]_x
V^{(m)}, V^{(m)}
)dx
\\
&\equiv 
-\|
\left(
g(u_x,u_x)
\right)^{1/2}
\nabla_xV^{(m)}
\|_{L^2}^2.
\end{aligned}
\label{equation:E7}
\end{equation}
As for the term 
$m\left(\nabla_xJ\right)\nabla_xV^{(m)}$,
note first that there exists a positive constant 
$C_1=C_1(N)>0$ such that 
\begin{equation}
\left|
\left(\nabla_xJ\right)
\right|(x)
\leqslant 
C_1(N)
\left(
g(u_x(x),u_x(x))
\right)^{1/2}
\label{equation:E8}
\end{equation}
holds uniformly with respect to $x$.
Thus we have 
\begin{equation}
\begin{aligned}
&\int_{\RR}
g(
m\left(\nabla_xJ\right)\nabla_xV^{(m)}, V^{(m)}
)dx
\\
&\leqslant
m
\|
\left(\nabla_xJ\right)\nabla_xV^{(m)}
\|_{L^2}
\|
V^{(m)}
\|_{L^2} 
\\
&\leqslant
mC_1(N)
\|
\left(
g(u_x,u_x)
\right)^{1/2}
\nabla_xV^{(m)}
\|_{L^2}
\|
V^{(m)}
\|_{L^2} 
\\
&
\leqslant
\rho
\|
\left(
g(u_x,u_x)
\right)^{1/2}
\nabla_xV^{(m)}
\|_{L^2}^2
+
\frac{m^2C_1^2}{4\rho}
\|
V^{(m)}
\|_{L^2}^2
\\
&\equiv
\rho
\|
\left(
g(u_x,u_x)
\right)^{1/2}
\nabla_xV^{(m)}
\|_{L^2}^2
\end{aligned}
\label{equation:E9}
\end{equation}
for any $\rho >0$. 
Note that the third inequality  above is 
due to the Schwartz inequality.
\par
In the same way, as for the term 
$-K_x\nabla_xJV^{(m)}$ and
$-K_xJ\nabla_xV^{(m)}$, we have
\begin{equation}
\begin{aligned}
&
\int_{\RR}
g( -K_x\nabla_xJV^{(m)}, V^{(m)} )dx
+
\int_{\RR}
g( -K_xJ\nabla_xV^{(m)}, V^{(m)} )dx
\\
&=
\int_{\RR}
g( K_{xx}JV^{(m)}, V^{(m)} )dx
+
2\int_{\RR}
g( K_xJV^{(m)}, \nabla_xV^{(m)} )dx
\\
&=
2\int_{\RR}
g( K_xJV^{(m)}, \nabla_xV^{(m)} )dx 
\\
&=
-\frac{2}{3a}
\int_{\RR}
g( 
\left(g(u_x,u_x)\right)^{1/2}JV^{(m)},
\left(g(u_x,u_x)\right)^{1/2}\nabla_xV^{(m)} 
)dx
\\
&
\leqslant
\rho
\|
\left(
g(u_x,u_x)
\right)^{1/2}
\nabla_xV^{(m)}
\|_{L^2}^2
+
\left(\frac{2}{3|a|}\right)^2\frac{1}{4\rho}
\|
\left(
g(u_x,u_x)
\right)^{1/2}JV^{(m)}
\|_{L^2}^2
\\
&\equiv
\rho
\|
\left(
g(u_x,u_x)
\right)^{1/2}
\nabla_xV^{(m)}
\|_{L^2}^2
\end{aligned}
\label{equation:E10}
\end{equation}
for any $\rho >0$.
\par
By combining 
\eqref{equation:E3},
\eqref{equation:E4},
\eqref{equation:E5},
\eqref{equation:E6},
\eqref{equation:E7},
\eqref{equation:E9}
and
\eqref{equation:E10}, 
and by taking  $\rho =1/4$,
we deduce 
\begin{equation}
\int_{\RR}
g
(
F_2(t), V^{(m)}(t)
)dx
\equiv
-\frac{1}{2}
\|
\left(
g(u_x(t),u_x(t))
\right)^{1/2}
\nabla_xV^{(m)}(t)
\|_{L^2}^2.
\label{equation:E11}
\end{equation}
\par
Thirdly, we consider $F_3$.
There never appear the terms containing  higher ordered derivative 
like $\nabla_x^{m+l}u_x$ with $l\in \mathbb{N}$ in $F_3$.
Hence it is easy to obtain that
\begin{equation}
\begin{aligned}
\int_{\RR}
g(
F_3(t),
V^{(m)}(t)
)dx  
&
\leqslant 
C(a,b,m,N,r_0)
P(N^{\ep}_4(t)+N^{\ep}_{m-1}(t))
N^{\ep}_{m}(t)
\|
V^{(m)}(t)
\|_{L^2} 
\\
&\equiv 0. 
\end{aligned}
\label{equation:E2}
\end{equation}
Here we add some comments on the estimation. 
The curvature tensor 
is estimated as follows: 
for $l\geqslant 0$ (resp. $j\geqslant 1$ ) 
and $U,V,W\in \Gamma(u^{-1}TN)$, 
there exists a positive constant 
$C(N,l)>0$ (resp. $C(N,j)>0$ ) 
such that
\begin{align}
\left|
\nabla_x^l
\left[ R(U,V) W\right]
\right|(x)
&\leqslant 
C(N,l)
\sum_{\begin{smallmatrix}
       p+q+r+j=l \\
       p,q,r,j \geqslant 0
       \end{smallmatrix}}
\left| (\nabla_x^jR)  \right|
\left| \nabla_x^pU   \right|
\left| \nabla_x^qV   \right|
\left| \nabla_x^rW   \right|(x),
\nonumber 
\\
\left| (\nabla_x^jR)  \right|(x)
&\leqslant 
C(N,j)
\sum_{\alpha =1}^j
\sum_{\begin{smallmatrix}
       \alpha + \sum_{h=1}^{\alpha}p_h=j \\
       p_h \geqslant 0
       \end{smallmatrix}}
\left| \nabla_x^{p_1}u_x  \right|
\cdots
\left| \nabla_x^{p_{\alpha}}u_x \right|(x)
\label{equation:EStr}
\intertext{uniformly with respect to $x$, 
where $\left|\cdot\right|=\left(g(\cdot,\cdot)\right)^{1/2}$.
Similarly, 
the $(1,1)$-tensor field 
$(\nabla_x^{j+1}J)$ with $j\geqslant 0$ is estimated as }
\left| (\nabla_x^{j+1}J)  \right|(x)
&\leqslant 
C(N,j)
\sum_{\alpha =1}^{j+1}
\sum_{\begin{smallmatrix}
       \alpha + \sum_{h=1}^{\alpha}p_h=j+1 \\
       p_h \geqslant 0
       \end{smallmatrix}}
\left| \nabla_x^{p_1}u_x  \right|
\cdots
\left| \nabla_x^{p_{\alpha}}u_x \right|(x)
\end{align}
for some positive constant $C(N,j)>0$.
Observing them, 
we can see that higher ordered derivatives 
never appear in $F_3$ and thus \eqref{equation:E2} 
is obtained.
Note also $K_tV^{(m)}$ is contained in $F_3$. 
The requirement $m\geqslant 4$  
comes to control this term. 
In other words, the $L^{\infty}$-norm of 
$K_t$ is bounded by some positive constant 
$C=C(a,r_0)$. Hence 
$K_tV^{(m)}$ is also harmless 
in the 
estimation \eqref{equation:E2}.
\par
Finally we consider the term $\ep F_1$. 
By repeatedly using integration by parts and 
the Schwartz inequality as before, 
it is easy to check that 
\begin{equation}
\int_{\RR}
g(
\ep F_1(t), V^{(m)}(t)
)dx
\equiv
\rho
\ep
\|
\nabla_x^2V^{(m)}(t)
\|_{L^2}^2
\label{equation:E12}
\end{equation}
for any $\rho >0$.
Thus, by taking $\rho =1/2$, it follows from 
\eqref{equation:EE1}
and 
\eqref{equation:E12}
that 
\begin{equation}
\int_{\RR}
g(
-\ep \nabla_x^4V^{(m)}(t)+\ep F_1(t), 
V^{(m)}(t)
)dx
\equiv
-\frac{\ep}{2}
\|
\nabla_x^2V^{(m)}(t)
\|_{L^2}^2.
\label{equation:E13}
\end{equation}
Consequently, 
\eqref{equation:EE2},
\eqref{equation:EE3},
\eqref{equation:E11},
\eqref{equation:E2},
and 
\eqref{equation:E13}
yield that 
\eqref{equation:E1} 
is estimated as follows: 
\begin{equation}
\begin{aligned}
&\frac{1}{2}\frac{d}{dt}
\|
V^{(m)}(t)
\|_{L^2}^2
+
\frac{\ep}{2}
\|\nabla_x^2 V^{(m)}(t) \|_{L^2}^2  
+
\frac{1}{2}
\|\left( g(u_x(t),u_x(t))\right)^{1/2}
\nabla_xV^{ (m)} (t) \|_{L^2}^2 
\\
&\leqslant 
C(a,b,m,N,r_0)
P(N^{\ep}_4(t)+N^{\ep}_{m-1}(t))
\left(
N^{\ep}_m(t)
\right)^2
\end{aligned}
\label{equation:E14}
\end{equation}
for some 
$C(a,b,m,N,r_0)>0$
and increasing function $P(\cdot)$.
\par 
On the other hands, 
it is easy to prove 
\begin{align}
&\frac{1}{2}\frac{d}{dt}
\left\|
u_x(t)
\right\|_{H^{m-1}}^2
+
\frac{\ep}{2}
\sum_{l=0}^{m-1}
\|\nabla_x^{l+2}u_x(t) \|_{L^2}^2
\equiv 
0.
\label{equation:E15}
\end{align}
By adding \eqref{equation:E14} and \eqref{equation:E15}, 
we obtain the desired estimate 
\eqref{equation:Energyinequality}.
\qed
\end{bew}
Lemma~\ref{lemma:Energy} follows immediately 
from Proposition~\ref{proposition:Energyestimate}
in the following way. 
If $m=4$, 
then \eqref{equation:Energyinequality} implies that 
$$
\left(
N^{\ep}_4(t)
\right)^2
\leqslant 
r_0^2
\exp
\left( 
2C(a,b,4,N,r_0)t
\right)
\quad 
\text{for}
\quad
t\in [0,T_{\ep}^{*}].
$$
If we set $t=T_{\ep}^{*}$, then this becomes 
$$
4r_0^2
=
\left(
N^{\ep}_4(T_{\ep}^{*})
\right)^2
\leqslant
r_0^2
\exp
\left( 
2C(a,b,4,N,r_0)T_{\ep}^{*}
\right),
$$
which implies 
$$
T_{\ep}^{*}
\geqslant
T
\equiv
\frac{2C(a,b,4,N,r_0)}{\log 4}.
$$
Clearly $T$ depends only on $a,b,N, \|u_{0x}\|_{H^4}$, 
being independent of $\ep\in (0, 1)$,
and 
$\{u^{\ep}_x\}_{\ep\in (0, 1)}$
is a bounded sequence in
$L^{\infty}(0,T;H^4(\RR;TN))$.
Then, by using the Gronwall inequality for 
$m=5,6,\ldots$ inductively, 
we obtain that 
$\{u^{\ep}_x\}_{\ep\in (0, 1)}$
is a bounded sequence in
$L^{\infty}(0,T;H^m(\RR;TN))$.
\qed
\end{bew}
\begin{remark}
\label{remark:remark3}
$\{u_x^{\ep}\}_{\ep \in (0,1)}$ 
gains the regularity in the following sense: 
By integrating \eqref{equation:Energyinequality} 
on $[0, T]$,  
we obtain
$$
\frac{\ep}{2}
\left(
\|\nabla_x^2 V^{\ep, (m)} \|_{L^2((0,T)\times\RR;TN)}^2
+
\sum_{l=0}^{m-1}
\|\nabla_x^{l+2}u_x^{\ep} \|_{L^2((0,T)\times\RR;TN)}^2
\right)
\leqslant 
C
$$
for some constant $C=C(a, b, N, \|u_{0x}\|_{H^m},T)>0$ 
independent of $\ep \in (0,1)$. 
This implies that the sequence
$
\{\ep^{1/2}
\nabla_x^mu_x^{\ep}\}_{\ep \in (0,1)}
$
is bounded in 
$L^2(0,T;H^2(\RR;TN))$. 
From this and Lemma~\ref{lemma:Energy} 
it is obvious that
$
\{u_t^{\ep}\}_{\ep \in (0,1)}
$
is also a bounded sequence in
$L^2(0,T;H^{m-2}(\RR;TN))$. 
We will use this property in 
the compactness argument 
in the next section.
\end{remark}
\section{Proof of Theorem \ref{theorem:eo}}
\label{section:proof1}
\begin{bew}{of Theorem~\ref{theorem:eo}}
We are now in a position to complete the proof of 
Theorem~\ref{theorem:eo}. 
We have only to solve 
\eqref{equation:pde}-\eqref{equation:data} 
in the positive direction of the time variable.
\begin{bew}{of existence}
Suppose that $u_0\in H^{m+1}(\RR;N)$  
with the integer $m\geqslant 4$ is given. 
By applying 
Proposition~\ref{proposition:pr}  
as $k=m$,
we construct a sequence  
$\{u^{\ep}\}_{\ep \in (0,1)}$ 
solving \eqref{equation:pde4}-\eqref{equation:data4} 
for each $\ep>0$. 
Recall that  Lemma~\ref{lemma:Energy} implies that
there exists  
$T=T(a, b, N,\|u_{0x}\|_{H^4(\RR;TN)})>0$  
which is independent of $\ep \in (0,1)$  
such that 
$\{u^{\ep}_x\}_{\ep \in (0,1)}$ 
is bounded in 
$L^{\infty}(0,T;H^m(\RR;TN))$. 
Recall also, as stated in 
Remark~\ref{remark:remark3} 
in the  previous section,
$
\{u_t^{\ep}\}_{\ep \in (0,1)}
$
is bounded in the class
$L^2(0,T;H^{m-2}(\RR;TN))$. 
Having them in mind, define $v^{\ep}=w{\circ}u^{\ep}$.
Then the boundnesses above imply respectively 
that
$\{v^{\ep}_x\}_{\ep \in (0,1)}$ 
is bounded in 
$L^{\infty}(0,T;H^m(\RR;\RR^d))$ 
and 
$\{v^{\ep}_t\}_{\ep \in (0,1)}$ 
is bounded in 
$L^{2}(0,T;H^{m-2}(\RR;\RR^d))$.   
Especially,  this boundness of 
$\{v^{\ep}_t\}_{\ep \in (0,1)}$ 
yields that 
$\{v^{\ep}_x\}_{\ep \in (0,1)}$ 
is bounded in the class
$C^{0,1/2}([0,T];H^{m-3}(\RR;\RR^d))$.
Then the standard compactness arguments imply 
that there exists a subsequence 
$\{v^{j}\}_{j\in \mathbb{N}}$ 
and $v$ such that 
\begin{alignat}{4}
& v^j_x
\stackrel{w^{\star}}{\longrightarrow}
v_x
\quad
&
\text{in}
\quad
&
L^{\infty}(0,T;H^{m}(\RR;\RR^d))
\quad
&
\text{as}
\quad
&
j\to \infty,
\label{equation:converge1}
\\
& v^j_x
\longrightarrow
v_x
\quad
&
\text{in}
\quad
&
C([0,T];H^{m-1}_{loc}(\RR;\RR^d))
\quad
&
\text{as}
\quad
&
j\to \infty,
\label{equation:converge2} 
\\
& v^j
\longrightarrow
v
\quad
&
\text{in}
\quad
&
C([0,T]\times \overline{B(0,R)};\RR^d))
\quad
&
\text{as}
\quad
&
j\to \infty 
\label{equation:converge3}
\end{alignat}
for any $R>0$, 
where 
$\overline{B(0,R)}=\left\{ x\in \RR \ | \ |x|\leqslant R  \right\}$.
In particular, 
\eqref{equation:converge3} implies that 
$v\in C([0,T]\times \RR;w(N))$ 
and  $w^{-1}{\circ}v$ satisfies the initial condition 
\eqref{equation:data}.
Furthermore, it is easy to check that $v$ satisfies 
\eqref{equation:pde5} with $\ep =0$. 
At this time, notice that
$v_x\in 
L^{\infty}(0,T;H^{m}(\RR;\RR^d))
\cap
C([0,T];H^{m-1}(\RR;\RR^d))
$
follows.
As a consequence, we have   
$u=w^{-1}{\circ} v\in C([0,T]\times \RR;N)$ 
with 
\begin{equation}
\label{equation:exsol}
u_x\in
L^{\infty}(0,T;H^{m}(\RR;TN))
\cap
C([0,T];H^{m-1}(\RR;TN))
\end{equation}
which solves \eqref{equation:pde} 
with the initial data $u_0$. 
Thus we  complete the proof of the existence
of time-local solutions. 
\qed
\end{bew}
\begin{remark}
\label{remark:L2}
For the solution $u=w^{-1}{\circ} v$, 
since
$v_x\in L^{\infty}(0,T;H^{m}(\RR;\RR^d))$, 
$v_t$ belongs to  $L^{\infty}(0,T;H^{m-2}(\RR;\RR^d))$, 
and thus
we see that $v-w{\circ}u_0$
belongs to 
$C^{0,1}([0,T];H^{m-2}(\RR;\RR^d))$. 
\end{remark}
\begin{bew}{of uniqueness}
Let $u,v\in C([0,T]\times \RR;N)$ 
be solutions of \eqref{equation:pde}-\eqref{equation:data} 
with \eqref{equation:exsol}, 
and let $u(0,x)=v(0,x)$. 
Identify $u,v$ with $w{\circ} u, w{\circ} v$. 
Then $u$ and $v$ satisfy 
$$v_t-av_{xxx}=f(v, v_x, v_{xx}),$$
where 
\begin{align}
f(v, v_x, v_{xx})
=&
a\left\{
\left[A(v)(v_x,v_x)\right]_x
+
A(v)(v_{xx}+A(v)(v_x,v_x), v_x)
\right\}
\nonumber \\
&+
dw_{w^{-1}\circ v}J_{w^{-1}\circ v}dw^{-1}_v(v_{xx}+A(v)(v_x,v_x))
+
b\left|v_x\right|^2v_x
\nonumber
\end{align}
for $v:\RR\to N$. 
As is stated in Remark~\ref{remark:L2},
both $u-w{\circ}u_0$ and $ v-w{\circ}u_0$ 
belong to the class 
$C^{0,1}([0,T];H^{m-2}(\RR;\RR^d))$ 
and thus 
$z=u-v$ is well-defined as a $\RR^d$-valued function. 
Taking the difference between two equations, 
we have 
$$
z_t-az_{xxx}=f(u, u_x, u_{xx})-f(v, v_x, v_{xx}), 
$$
To prove that $z=0$, we can show that there exists 
a constant $C>0$  
depending only on  
$a, b, N$, $\|u_x\|_{L^{\infty}(0,T;H^2(\RR;\RR^d))}$, and
$\|v_x\|_{L^{\infty}(0,T;H^2(\RR;\RR^d))}$
such that
\begin{equation}
\frac{d}{dt}
\|z(t)\|_{H^1(\RR;\RR^d)}^2
\leqslant
C\|z(t)\|_{H^1(\RR;\RR^d)}^2.
\label{equation:deff}
\end{equation}
This estimate can be 
obtained by  
completely same calculation
as that in the proof of the uniqueness in \cite{Onodera1}.
Note, though the only case that 
$(N,J,g)$ is a K\"ahler manifold is discussed 
in \cite{Onodera1}, 
the argument proving the uniqueness works 
also when  $(N,J,g)$ is a compact 
almost Hermitian manifold.
Thus we omit the proof of \eqref{equation:deff}.
\qed
\end{bew}
\begin{bew}{of the continuity in time of \ $\nabla_x^mu_x$
in $L^2(\RR;TN)$}
We have already proved the existence of a unique solution 
$u\in C([0,T]\times \RR;N)$ 
with \eqref{equation:exsol}. 
Thus the proof of 
$\nabla_x^mu_x\in C([0,T];L^2(\RR;TN)$ 
is left.
Let $v=w{\circ} u$. 
To obtain this continuity, 
it suffices to show that 
$dw_u(V^{(m)})$ belongs to 
 $C([0,T];L^2(\RR;\RR^d))$.
\par 
First of all, 
the energy estimate 
\eqref{equation:Energyinequality} implies
$\left( d/dt \right)
\left( N^{\ep}_m(t)  \right)^2
\leqslant C$
for some $C>0$ which is independent of $\ep\in (0,1)$. 
Hence we deduce 
\begin{align}
&\|V^{\ep, (m)}(t)\|_{L^2(\RR;TN)}^2
+
\|u_x^{\ep}(t)\|_{H^{m-1}(\RR;TN)}^2 
\nonumber 
\\
&\qquad \qquad \leqslant
\|V^{\ep, (m)}(0)\|_{L^2(\RR;TN)}^2
+
\|u_x^{\ep}(0)\|_{H^{m-1}(\RR;TN)}^2
+Ct.
\nonumber 
\intertext{
Letting $\ep \downarrow 0$, we see that
$V^{(m)}(t)=(e^{K}\nabla_x^mu_x)(t)\in L^2(\RR;\RR^d)$ 
makes sense for all $t\in [0, T]$, and }
&\|V^{ (m)}(t)\|_{L^2(\RR;TN)}^2
+
\|u_x(t)\|_{H^{m-1}(\RR;TN)}^2 
\nonumber 
\\
&\qquad \qquad \leqslant
\|V^{ (m)}(0)\|_{L^2(\RR;TN)}^2
+
\|u_x(0)\|_{H^{m-1}(\RR;TN)}^2
+Ct.
\nonumber
\end{align}
Noting that $u_x\in C([0,T];H^{m-1}(\RR;TN))$, we have 
\begin{align}
\limsup_{t\to 0}
\| V^{ (m)} (t)\|_{L^2(\RR;TN)}^2
&\leqslant
\|V^{ (m)}(0)\|_{L^2(\RR;TN)}^2.
\label{equation:es11}
\intertext{Since $w$ is the isometric embedding, 
\eqref{equation:es11} is equivalent to}
\limsup_{t\to 0}
\| dw_u ( V^{ (m)} ) (t)\|_{L^2(\RR;\RR^d)}^2
&\leqslant
\| dw_u ( V^{ (m)} ) (0)  \|_{L^2(\RR;\RR^d)}^2.
\label{equation:es1}
\intertext{Moreover, 
since $v_x\in L^{\infty}(0,T;H^{m}(\RR;\RR^d)) 
\cap 
C([0,T];H^{m-1}(\RR;\RR^d))$, 
we see 
$dw_u( V^{ (m)}  )(t)$ 
is weakly continuous in $L^2(\RR;\RR^d)$. 
Hence it follows that}
\|dw_u( V^{ (m)} )(0)\|_{L^2(\RR;\RR^d)}^2
&\leqslant
\liminf_{t\to 0}
\|dw_u( V^{ (m)}  )(t)\|_{L^2(\RR;\RR^d)}^2.
\label{equation:es2}
\intertext{From \eqref{equation:es1} and \eqref{equation:es2}, 
we obtain} 
\lim_{t\to 0}
\|dw_u( V^{ (m)}  )(t)\|_{L^2(\RR;\RR^d)}^2
&=
\|dw_u( V^{ (m)}   )(0)\|_{L^2(\RR;\RR^d)}^2.
\label{equation:es3}
\end{align}
Consequently, \eqref{equation:es3} 
and the weak continuity of $dw_u( V^{ (m)}  )(t)$ 
in the class $L^2(\RR;\RR^d)$ imply that
$dw_u( V^{ (m)}   )(t)$ 
is strongly continuous in $L^2(\RR;\RR^d)$ 
at $t=0$. 
By the uniqueness of $u$,  
we see 
$dw_u( V^{ (m)}  )(t)$ 
is strongly continuous at each $t\in [0, T]$ in the same way. 
Thus we complete the proof.
\qed
\end{bew}
\qed
\end{bew}
\section{Sketch of the 
proof of Theorem~\ref{theorem:eeo} and \ref{theorem:meo}}
\label{section:special}
This section is devoted to 
the outline of 
the proof of Theorem \ref{theorem:eeo} 
and \ref{theorem:meo}.
Recall in both cases, $N$ is 
supposed to be a compact K\"ahler manifold. 
\begin{bew}{of Theorem~\ref{theorem:eeo}}
Since $N$ is a compact K\"ahler manifold,
the procedures of the proof is almost parallel to that in 
\cite{Onodera1}. 
There is a difference to the proof of Theorem~\ref{theorem:eo} 
in the energy estimate. 
Due to the K\"ahler condition, 
the classical energy method works effectively. 
In other words, 
we do not need to use the gauge  transformation 
of $\nabla_x^mu_x$ used in the proof of Theorem~\ref{theorem:eo}.
This is the reason that 
this theorem holds for $m\geqslant 2$.
Indeed, we can obtain the following. 
\begin{lemma}
\label{lemma:classicalenergy}
Let $\{ u^{\ep}  \}_{\ep \in (0,1)}$ be a 
sequence of solution of 
\eqref{equation:pde4}-\eqref{equation:data4}
constructed in Proposition~\ref{proposition:pr} 
as $k=m\geqslant 2$.
Then there exists a constant  
$T>0$ 
depending only on 
$a, b, N$, and $\|u_{0x}\|_{H^2(\RR;TN)}$  
such that 
$\{u_x^{\ep}\}_{\ep\in (0, 1)}$ 
is bounded in 
$L^{\infty}(0, T; H^m(\RR;TN))$.
\end{lemma}
\begin{bew}{of Lemma~\ref{lemma:classicalenergy}}
By the completely same calculus as that 
in \cite[Lemma~4.1]{Onodera1}, 
we can show that 
\begin{align}
& \frac{d}{dt}
  \| u_{x}^{\ep}(t) \|_{H^2(\RR;TN)}^2
  \leqslant
  C(a,b,N)\sum_{r=4}^8\| u_x^{\ep}(t) \|_{H^2(\RR;TN)}^r,
\label{equation:energy2}
\\
& \frac{d}{dt}
  \| u_x^{\ep}(t) \|_{H^k(\RR;TN)}^2
  \leqslant 
  C(a,b,N, 
  \| u_x^{\ep}(t) \|_{H^{k-1}(\RR;TN)}
  )
  \| u_x^{\ep}(t) \|_{H^k(\RR;TN)}^2 
\label{equation:energyk}
\end{align}
for $3\leqslant k\leqslant m$ 
hold for all $t\in [0,T_{\ep}]$. 
From \eqref{equation:energy2} and  \eqref{equation:energyk}, 
the desired boundness is immediately obtained. 
See \cite[Lemma~4.1]{Onodera1} for details. 
\qed
\end{bew}
The other parts of the proof of Theorem~\ref{theorem:eeo} 
are same as that was discussed in Theorem~\ref{theorem:eo}.
Thus we omit the detail. 
\qed
\end{bew}
Next, let $(N, J, g)$ be a compact Riemann surface  
with constant Gaussian curvature $K$, 
and assume that $a\neq 0$ and $b=aK/2$. 
Theorem~\ref{theorem:eeo} tells us that, 
given a initial data 
$u_0\in H^{m+1}(\RR;N)$, 
there exists  
$T=T(a, b, N,\|u_{0x}\|_{H^2(\RR;TN)})>0$ 
such that the IVP 
\eqref{equation:pde}-\eqref{equation:data} 
admits a unique time-local solution  
$u\in C([0,T);H^{m+1}(\RR;N))$.
\par 
In what follows  
we will extend the existence time of $u$ 
over $[0,\infty)$. 
For this,  
we have the following energy conversation laws.
\begin{lemma}
\label{lemma:conserved}
For $u\in C([0,T);H^{m+1}(\RR;N))$ 
solving \eqref{equation:pde}-\eqref{equation:data}, 
the following quantities
\begin{align}
&\|u_x(t)\|_{L^2(\RR;TN)}^2,
\nonumber \\
&E(u(t))
=
\|\nabla_x^2u_x(t)\|_{L^2(\RR;TN)}^2
+
\frac{K^2}{8}
\int_{\RR}
\left(
g(u_x(t),u_x(t))
\right)^3dx
\nonumber \\
&\quad \phantom{E(u(t))=}-K
\int_{\RR}
\left(
g(u_x(t),\nabla_xu_x(t))
\right)^2dx
\nonumber \\
&\quad \quad \phantom{E(u(t))=}-\frac{3K}{2}
\int_{\RR}
g(u_x(t),u_x(t))
g(\nabla_xu_x(t),\nabla_xu_x(t))dx
\nonumber
\end{align}
are preserved with respect to $t\in [0,T)$.
\end{lemma}
\begin{bew}{of Lemma~\ref{lemma:conserved}} 
The proof is also  same as that was discussed in 
\cite[Lemma~6.1]{Onodera1}. 
Thus we omit the detail.
\qed
\end{bew}
\begin{bew}{of Theorem~\ref{theorem:meo}}
Let 
$u\in C([0,T);H^{m+1}(\RR;N))$ 
be a time-local solution of  
\eqref{equation:pde}-\eqref{equation:data} 
which exists on the maximal time interval $[0,T)$. 
If $T=\infty$, 
Theorem~\ref{theorem:meo} holds true. 
Thus we only need to consider the case $T<\infty$. 
From Lemma~\ref{lemma:conserved}, 
we know that
\begin{equation}
\|u_x(t)\|_{L^2(\RR;TN)}^2
=
\|u_{0x}\|_{L^2(\RR;TN)}^2,
\quad
E(u(t))=E(u_0).
\label{equation:q0}
\end{equation}
Hence it follows that
\begin{align}
\|\nabla_x^2u_x(t)\|_{L^2(\RR;TN)}^2
=&
E(u_0)-\frac{K^2}{8}
\int_{\RR}
\left(
g(u_x(t),u_x(t))
\right)^3dx
\nonumber \\
&\quad +
K
\int_{\RR}
\left(
g(u_x(t),\nabla_xu_x(t))
\right)^2dx
\nonumber \\
&\quad \quad +
\frac{3K}{2}
\int_{\RR}
g(u_x(t),u_x(t))
g(\nabla_xu_x(t),\nabla_xu_x(t))dx
\nonumber \\
\leqslant&
E(u_0)
+
C|K|
\|u_x(t)\|_{L^{\infty}(\RR;TN)}^2
\|\nabla_xu_x(t)\|_{L^2(\RR;TN)}^2.
\nonumber 
\end{align}
The second term of the right hand side of the above is 
estimated as follows. 
At first, we have
\begin{align}
\|\nabla_xu_x(t)\|_{L^2(\RR;TN)}^2 
&= 
-\int_{\RR}g\left( u_x(t), \nabla_x^2u_x(t)  \right)dx
\nonumber 
\\
&\leqslant
\|u_x(t)\|_{L^2(\RR;TN)}
\|\nabla_x^2u_x(t)\|_{L^2(\RR;TN)}
\nonumber 
\\
&=
\|u_{0x}\|_{L^2(\RR;TN)}
\|\nabla_x^2u_x(t)\|_{L^2(\RR;TN)}.
\label{equation:gn}
\end{align}
Next, note that 
$dw_u(\nabla_xu_x)=v_{xx}+A(v)(v_x,v_x)$ 
holds for $v=w{\circ}u$ 
by the definition of the covariant derivative 
along the mapping $u$.
By noting this
 and by using 
\eqref{equation:gn} and Sobolev's inequality,
we obtain 
\begin{align}
&
\|u_x(t)\|_{L^{\infty}(\RR;TN)}^2 
\nonumber 
\\
&=
\|v_x(t)\|_{L^{\infty}(\RR;\RR^d)}^2 
\nonumber 
\\
&\leqslant
C
\|v_x(t)\|_{L^2(\RR;\RR^d)}
\|v_{xx}(t)\|_{L^2(\RR;\RR^d)}
\nonumber 
\\
&\leqslant
C
\|v_x(t)\|_{L^2(\RR;\RR^d)} 
\nonumber 
\\
&\qquad \times
\bigg( 
\| v_{xx}(t)+A(v)(v_x,v_x)(t) \|_{L^2(\RR;\RR^d)}
+
\| A(v)(v_x,v_x)(t) \|_{L^2(\RR;\RR^d)}
\bigg)
\nonumber 
\\
&\leqslant
C
\|v_x(t)\|_{L^2(\RR;\RR^d)}
\nonumber 
\\
&\qquad \times
\bigg( 
\| v_{xx}(t)+A(v)(v_x,v_x)(t) \|_{L^2(\RR;\RR^d)} 
\nonumber 
\\
&\qquad \qquad \qquad \qquad \qquad \qquad
+
C(N)
\| v_x(t) \|_{L^{\infty}(\RR;\RR^d)}
\| v_x(t) \|_{L^2(\RR;\RR^d)}
\bigg)
\nonumber 
\\
&=
C
\|u_x(t)\|_{L^2(\RR;TN)} 
\nonumber 
\\
&\qquad \times
\bigg( 
\| \nabla_xu_x(t) \|_{L^2(\RR;TN)}
+
C(N)
\| u_x(t) \|_{L^{\infty}(\RR;TN)}
\| u_x(t) \|_{L^2(\RR;TN)}
\bigg)
\nonumber 
\\
&\leqslant 
C
\|u_{0x}\|_{L^2(\RR;TN)} 
\nonumber 
\\
&\qquad \times
\bigg( 
\| u_{0x} \|_{L^2(\RR;TN)}^{1/2}
\| \nabla_x^2u_x(t) \|_{L^2(\RR;TN)}^{1/2} 
\nonumber 
\\
& \qquad \qquad \qquad \qquad \qquad \qquad
+
C(N)
\| u_x(t) \|_{L^{\infty}(\RR;TN)}
\| u_{0x} \|_{L^2(\RR;TN)}
\bigg)
\nonumber 
\\
&= 
C
\| u_{0x} \|_{L^2(\RR;TN)}^{3/2}
\| \nabla_x^2u_x(t) \|_{L^2(\RR;TN)}^{1/2} 
\nonumber 
\\
& \qquad \qquad \qquad +
C(N)
\| u_x(t) \|_{L^{\infty}(\RR;TN)}
\| u_{0x} \|_{L^2(\RR;TN)}^2,
\nonumber
\end{align}
which implies 
\begin{equation} 
\| u_x(t) \|_{L^{\infty}(\RR;TN)}
\leqslant 
C(N,\| u_{0x} \|_{L^2(\RR;TN)} )
\left( 
1+\| \nabla_x^2u_x(t) \|_{L^2(\RR;TN)}^{1/4}
\right).
\label{equation:so}
\end{equation}
From \eqref{equation:q0}, \eqref{equation:gn} and \eqref{equation:so},  
we deduce
\begin{align}
&\|\nabla_x^2u_x(t)\|_{L^2(\RR;TN)}^2
\nonumber 
\\
& \leqslant
E(u_0)
+
C(K,N,\|u_{0x}\|_{L^2(\RR;TN)}) 
\nonumber  
\\
&\qquad \qquad \times 
\left( 1+\|\nabla_x^2u_x(t)\|_{L^2(\RR;TN)}^{1/2}\right)
\|\nabla_x^2u_x(t)\|_{L^2(\RR;TN)}.
\nonumber
\end{align}
Thus $X=X(t)=1+\|\nabla_x^2u_x(t)\|_{L^2(\RR;TN)}^2$ 
satisfies 
$$X\leqslant 1+E(u_0)+C(K,N,\|u_{0x}\|_{L^2(\RR;TN)})X^{3/4},$$  
which implies that 
$X(t)$ is bounded, and thus 
\begin{align}
\sup_{t\in [0, T)}\|\nabla_x^2u_x(t)\|_{L^2(\RR;TN)}
&\leqslant
C(K, N,\|u_{0x}\|_{H^2(\RR;TN)})
\label{equation:q1}
\intertext{ for some $C=C(K, N,\|u_{0x}\|_{H^2(\RR;TN)})>0$. 
Interpolating \eqref{equation:q0} and \eqref{equation:q1} 
we have}
\sup_{t\in [0,T)}
\|u_x(t)\|_{H^2(\RR;TN)}
&\leqslant
C(K, N,\|u_{0x}\|_{H^2(\RR;TN)}).
\nonumber
\end{align}
Once we obtain the $H^2(\RR;TN)$-boundness 
of $u_x$, 
the desired $H^m(\RR;TN)$-boundness 
of $u_x$
follows from the use of \eqref{equation:energyk} 
inductively with respect to $k=3, \ldots, m$. 
Thus the existence time of $u$ can be extended beyond $T$. 
\qed
\end{bew}
\section{Appendix}
\label{section:appendix}
In this section, we 
check \eqref{equation:Equation} 
used in Section~\ref{section:energy}. 
Operating 
$e^K\nabla_x^{m+1}$ 
on the equation \eqref{equation:pde4}, 
we have 
\begin{align}
e^K\nabla_x^{m+1}u_t 
&=
-\ep\, e^K\nabla_x^{m+4}u_x
+a\, e^K\nabla_x^{m+3}u_x 
+e^K\nabla_x^{m+1}J\nabla_xu_x 
\nonumber 
\\
&
\qquad \qquad \qquad \qquad \qquad \qquad
+b\, e^K\nabla_x^{m+1}g(u_x,u_x)u_x.
\label{equation:AP0}
\end{align}
\par
First, to compute each term of \eqref{equation:AP0}, 
we use the following relation 
\begin{align}
e^K\nabla_x^{m+k}u_x
&=
\nabla_x
\left(e^K\nabla_x^{m+k-1}u_x \right)
-
K_xe^K\nabla_x^{m+k-1}u_x
\quad
\text{for}
\quad
k\in \mathbb{N}.
\label{equation:APx}
\end{align}
By using this relation repeatedly, we deduce
\begin{align}
e^K\nabla_x^{m+1}u_x
=&
\nabla_xV^{(m)}-K_xV^{(m)},
\label{equation:AP1}
\\
e^K\nabla_x^{m+2}u_x
=&
\nabla_x^2V^{(m)}
-2K_x\nabla_xV^{(m)}
-\left(
K_{xx}-K_x^2
\right)V^{(m)},
\label{equation:AP2}
\\
e^K\nabla_x^{m+3}u_x
=&
\nabla_x^3V^{(m)}
-3K_x\nabla_x^2V^{(m)}
-
3\left(
K_{xx}-K_x^2
\right)\nabla_xV^{(m)}
\nonumber 
\\
&-
\left(
K_{xxx}-3K_xK_{xx}+K_x^3
\right)V^{(m)},
\label{equation:AP3}
\\
e^K\nabla_x^{m+4}u_x
=&
\nabla_x^4V^{(m)}
-
4K_x\nabla_x^3V^{(m)}
-
6\left(
K_{xx}-K_x^2
\right)\nabla_x^2V^{(m)}
\nonumber 
\\
&-
4\left(
K_{xxx}-3K_xK_{xx}+K_x^3
\right)\nabla_xV^{(m)}
\nonumber 
\\
&
-
\left(
K_{xxxx}-4K_xK_{xxx}-3K_{xx}^2
+6K_x^2K_{xx}
-K_{x}^4
\right)V^{(m)}.
\label{equation:AP4}
\end{align}
\par
Moreover, \eqref{equation:AP1} and the 
Leibniz rule yield that 
\begin{equation}
\begin{aligned}
&e^K\nabla_x^{m+1}
\left[
g(u_x,u_x)u_x
\right] 
\\
&=     
2e^Kg(\nabla_x^{m+1}u_x, u_x)u_x
+
e^Kg(u_x, u_x)\nabla_x^{m+1}u_x 
\\
&\quad +
         \sum_{\begin{smallmatrix}
        \alpha+\beta+\gamma=m+1 \\
        \alpha, \beta, \gamma\geqslant 0\\
        \max\{\alpha, \beta, \gamma\}\leqslant m
       \end{smallmatrix}}
  \frac{(m+1)!}{\alpha!\beta!\gamma!}
  e^Kg(
  \nabla_x^{\alpha}u_x, \nabla_x^{\beta}u_x
  )
  \nabla_x^{\gamma}u_x
\\
&=
2g(\nabla_xV^{(m)}, u_x)u_x
+
g(u_x, u_x)\nabla_xV^{(m)} 
\\
&\quad \quad -
2g(K_xV^{(m)}, u_x)u_x
-
g(u_x, u_x)K_xV^{(m)} 
\\
&\quad \quad \quad +
         \sum_{\begin{smallmatrix}
        \alpha+\beta+\gamma=m+1 \\
        \alpha, \beta, \gamma\geqslant 0\\
        \max\{\alpha, \beta, \gamma\}\leqslant m
       \end{smallmatrix}}
  \frac{(m+1)!}{\alpha!\beta!\gamma!}
  e^Kg(
  \nabla_x^{\alpha}u_x, \nabla_x^{\beta}u_x
  )
  \nabla_x^{\gamma}u_x.
\end{aligned}
\label{equation:APb}
\end{equation}
\par
Next, we compute 
$e^K\nabla_x^{m+1}u_t$.
Note that
$$
\nabla_tu_x=\nabla_xu_t
\quad 
\text{and} 
\quad
\nabla_x\nabla_tu_x=\nabla_t\nabla_xu_x+R(u_x,u_t)u_x 
$$
follow from the definition of the Levi-Civita connection.
Using these commutative relations inductively, 
we have
\begin{equation}
\nabla_x^{m+1}u_t
=
\nabla_t\nabla_x^mu_x
+
\sum_{l=0}^{m-1}
\nabla_x^l
\left[ 
R(u_x, u_t)\nabla_x^{m-(l+1)}u_x 
\right].
\label{equation:APt1}
\end{equation}
By multiplying $e^K$ with 
\eqref{equation:APt1}, we have 
\begin{equation}
e^K\nabla_x^{m+1}u_t
=
e^K\nabla_t\nabla_x^mu_x
+
\sum_{l=0}^{m-1}
e^K\nabla_x^l
\left[ 
R(u_x, u_t)\nabla_x^{m-(l+1)}u_x 
\right].
\label{equation:APt2}
\end{equation}
By noting 
$e^K\nabla_t\nabla_x^mu_x
=\nabla_t\left(
e^K\nabla_x^mu_x
\right)-K_t\nabla_x^mu_x
=\nabla_tV^{(m)}-K_tV^{(m)}$, 
and by
substituting \eqref{equation:pde4} 
into the second term of \eqref{equation:APt2}, 
we deduce 
\begin{equation}
\begin{aligned}
&e^K\nabla_x^{m+1}u_t
=
\nabla_tV^{(m)}-K_tV^{(m)}
-\ep 
\sum_{l=0}^{m-1}
e^K\nabla_x^l
\left[ 
R(u_x, \nabla_x^3u_x)\nabla_x^{m-(l+1)}u_x 
\right]
\\
& \qquad \qquad \qquad  \quad +
a\sum_{l=0}^{m-1}
e^K\nabla_x^l
\left[ 
R(u_x, \nabla_x^2u_x)\nabla_x^{m-(l+1)}u_x 
\right]
\\
&  \qquad \qquad \qquad \quad \quad +
\sum_{l=0}^{m-1}
e^K\nabla_x^l
\left[ 
R(u_x, J\nabla_xu_x)\nabla_x^{m-(l+1)}u_x 
\right].
\end{aligned}
\label{equation:APt3}
\end{equation}
(Note that 
$R(u_x, b\,g(u_x,u_x)u_x)\nabla_x^{m-(l+1)}u_x=0$ 
since $R(u_x, u_x)=0$.)
The fourth term of the right hand side of 
\eqref{equation:APt3} is 
decompose as 
\begin{equation}
\begin{aligned}
&
a\sum_{l=0}^{m-1}
e^K\nabla_x^l
\left[ 
R(u_x, \nabla_x^2u_x)\nabla_x^{m-(l+1)}u_x 
\right] 
\\
=&
a\left(
\sum_{l=0}^{m-1}
e^K\nabla_x^l
\left[
R(u_x, \nabla_x^2u_x)\nabla_x^{m-1-l}u_x
\right]
-R(u_x, \nabla_xV^{(m)})u_x
\right) 
\\
&+
a\, R(u_x, \nabla_xV^{(m)})u_x.
\end{aligned}
\label{equation:APt4}
\end{equation}
Note the term $\nabla_x^{m+1}u_x$ never appear 
in the first term of the right hand side of 
\eqref{equation:APt4}.
\par
Let us move to the computation of 
$e^K\nabla_x^{m+1}J\nabla_xu_x$. 
First, it follows from the definition that 
\begin{equation}
\left(\nabla_xJ\right)V
=
\nabla_xJV-J\nabla_xV
\quad
\text{for}
\quad
V\in \Gamma(u^{-1}TN),
\label{equation:tensorJ}
\end{equation}
where 
$\left(\nabla_xJ\right)$ 
is the covariant derivative of $(1,1)$-tensor $J$ 
with respect to $x$ along $u$ and is also
 $(1,1)$-tensor field along $u$. 
We will write 
$\left(\nabla_xJ\right)V$ 
not to be confused with 
$\nabla_xJ V$.
In the same way, 
$(\nabla_x^{j+1}J)$ with $j\geqslant 1$, 
which is the $(j+1)$-th covariant derivative 
of $(1,1)$-tensor field $J$,  
is also $(1,1)$-tensor field along $u$ 
defined inductively by the form
\begin{equation}
\left(\nabla_x^{j+1}J\right)V
=
\nabla_x\left(\nabla_x^{j}J\right) V
-\left(\nabla_x^{j}J\right)  \nabla_xV
\quad
\text{for}
\quad
V\in \Gamma(u^{-1}TN),
\nonumber
\end{equation}
where 
$\left(\nabla_x^{1}J\right)=\left(\nabla_xJ\right)$. 
Using \eqref{equation:tensorJ} repeatedly, 
we deduce 
\begin{align}
e^K\nabla_x^{m+1}J\nabla_xu_x
=
e^K\nabla_xJ\nabla_x^{m+1}u_x
+
e^K\sum_{l=1}^m
\nabla_x^l
\left(\nabla_xJ\right)
\nabla_x^{m+1-l}u_x
\label{equation:APJ1}
\end{align}
For the first term of the right hand side of 
\eqref{equation:APJ1}, 
\eqref{equation:AP1} and $e^KJ=Je^K$ 
yield
\begin{equation}
\begin{aligned} 
&e^K\nabla_xJ\nabla_x^{m+1}u_x 
\\
=&
\nabla_x
\left(
Je^K\nabla_x^{m+1}u_x
\right)
-K_x
Je^K\nabla_x^{m+1}u_x  
\\
=&
\nabla_xJ\nabla_xV^{(m)}
-K_x\nabla_xJV^{(m)}
-K_xJ\nabla_xV^{(m)}
-(K_{xx}-K_x^2)JV^{(m)}.
\end{aligned}
\label{equation:APJ2}
\end{equation}
For the second term of the right hand side of 
\eqref{equation:APJ1}, 
by regarding 
$\left(\nabla_xJ\right)$ and  $\nabla_x^{m+1-l}u_x$
as
a $(1,1)$-tensor field  and a $(1,0)$-tensor field respectively, 
we deduce
\begin{equation}
\begin{aligned}
&e^K\sum_{l=1}^m
\nabla_x^l
\left(\nabla_xJ\right)
\nabla_x^{m+1-l}u_x 
\\
&=
e^K\sum_{l=1}^m
\nabla_x^l 
C_1^2
\left(
\left(\nabla_xJ\right)
\otimes
\nabla_x^{m+1-l}u_x
\right) 
\\
&=
e^K\sum_{l=1}^m
C_1^2
\nabla_x^l 
\left(
\left(\nabla_xJ\right)
\otimes
\nabla_x^{m+1-l}u_x
\right) 
\\
&=
e^K\sum_{l=1}^m
C_1^2
\left\{
\sum_{j=0}^l
\frac{l!}{j!(l-j)!}
\left(\nabla_x^{j+1}J\right)
\otimes
\nabla_x^{m+1-l+(l-j)}u_x
\right\}
\\
&=
e^K\sum_{l=1}^m
\sum_{j=0}^l
\frac{l!}{j!(l-j)!}
\left(\nabla_x^{j+1}J\right)
\nabla_x^{m+1-j}u_x
\\
&=
m\, e^K\left(\nabla_xJ\right)\nabla_x^{m+1}u_x
+
e^K\sum_{l=1}^m
\sum_{j=1}^l
\frac{l!}{j!(l-j)!}
\left(\nabla_x^{j+1}J\right)
\nabla_x^{m+1-j}u_x,
\end{aligned}
\label{equation:APJ33}
\end{equation}
where 
$C_1^2:T_uN\otimes T_uN\otimes T_u^{*}N\to T_uN$ 
is a contraction which maps 
$x_i\otimes x_j\otimes y_k^{*}$ 
into 
$\sum_{j,k}y_k^{*}(x_j)x_i$.
Notice that the second equality of 
\eqref{equation:APJ33}
holds since the covariant derivative commutes with 
the contraction, 
and the third equality of \eqref{equation:APJ33}
is due to the fact that
$$\nabla_x
\left(
S\otimes T
\right)
=
(\nabla_xS)\otimes T+S\otimes (\nabla_xT)
$$
holds for any tensor $S$ and $T$. 
See, e.g., \cite{GHL} for these properties.
Moreover, by noting that
$f\left(\nabla_xJ\right)=\left(\nabla_xJ\right)f$ 
holds for any scalar function $f$ and by using \eqref{equation:AP1}, 
we deduce
\begin{equation}
\begin{aligned}
m\, e^K\left(\nabla_xJ\right)\nabla_x^{m+1}u_x
&=
m\, \left(\nabla_xJ\right)e^K\nabla_x^{m+1}u_x 
\\
&=
m\, \left(\nabla_xJ\right)\nabla_xV^{(m)}
-
mK_x\left(\nabla_xJ\right)V^{(m)}.
\end{aligned}
\label{equation:APJ4}
\end{equation}
Combining \eqref{equation:APJ1},\eqref{equation:APJ2},
\eqref{equation:APJ33},
and \eqref{equation:APJ4}, 
we obtain 
\begin{equation}
\begin{aligned}
&e^K\nabla_x^{m+1}J\nabla_xu_x 
\\
=&
\nabla_xJ\nabla_xV^{(m)}
-K_x\nabla_xJV^{(m)}-K_xJ\nabla_xV^{(m)}
-(K_{xx}-K_x^2)JV^{(m)}
\\
&\quad +
m\, \left(\nabla_xJ\right)\nabla_xV^{(m)}
-
mK_x\left(\nabla_xJ\right)V^{(m)} 
\\
&\quad \quad +
e^K\sum_{l=1}^m
\sum_{j=1}^l
\frac{l!}{j!(l-j)!}
\left(\nabla_x^{j+1}J\right)
\nabla_x^{m+1-j}u_x.
\end{aligned}
\label{equation:APJJ}
\end{equation}
Consequently, 
by substituting 
\eqref{equation:AP3},\eqref{equation:AP4},
\eqref{equation:APb},\eqref{equation:APt3},
\eqref{equation:APt4} and \eqref{equation:APJJ}
into 
\eqref{equation:AP0}, 
we deduce the desired equality  \eqref{equation:Equation}.
\\
\\
{\bf Acknowledgement.} \\
The author expresses gratitude to Hiroyuki Chihara  
for several discussions and valuable advice. 
Without them, this work  would not have been completed.

%
%
%
%

\begin{thebibliography}{10}
\bibitem{CSU}
Chang,~N.~H., Shatah,~J., Uhlenbeck,~K.:
{\it Schr\"odinger maps}. 
Comm.\  Pure Appl.\  Math.\  {\bf53} (2000), 590--602.

\bibitem{DR}
Da~Rios:
{\it On the motion of an unbounded fluid with a vortex filament 
     of any shape[in Italian]}.
 Rend, \ Circ.\ Mat.\ Palermo {\bf 22} (1906), 117--135. 

\bibitem{Ding}
Ding,~W.~Y.:
{\it On the Schr\"odinger flows}.
 Proceedings of the ICM,\ Vol.\ II.\ (2002), 283--291. 

\bibitem{Doi}Doi,~S.-I.: 
{\it Smoothing effects of Schr\"odinger evolution groups on 
Riemannian manifolds}. 
Duke Math.\ J. {\bf 82} (1996), 679--706.

\bibitem{FM}
Fukumoto,~Y., Miyazaki,~T.:
{\it Three-dimensional distortions of a	vortex filament with axial
	velocity}.
J.\ Fluid Mech.\ {\bf 222} (1991), 369--416.

\bibitem{GHL}
Gallot, ~S., Hulin, ~D., Lafontaine,~J.:
{\it Riemannian Geometry. }
third ed.,\ Universitext,\ Springer-Verlag,\ Berlin,\ 2004.

\bibitem{GR}
Gromov,~M.~L., Rohlin,~V.~A.:
{\it Embeddings and immersions in Riemannian geometry}.
 Usp.\ Mat.\ Nauk \ {\bf 25} (1970), 3--62(in\ Russian).
English \ translation:\ Russ.\ Math.\ Survey\
 {\bf 25} (1970), 1--57.   

\bibitem{Hasimoto}
Hasimoto,~H.:
{\it A soliton on a vortex filament}.
J.\ Fluid.\ Mech.\ {\bf 51} (51), 477--485.


\bibitem{Koiso}
Koiso,~N.:
{\it The vortex filament equation and a semilinear Schr\"odinger
	equation in a Hermitian symmetric space}.
 Osaka J.\ Math.\ {\bf 34} (1997), 199--214. 

\bibitem{LL}
Landau,~L.~D., Lifschitz,~E.~M.:
{\it On the theory of the dispersion of magnetic permeability in 
     ferro-aquatic bodies}.
Physica A (Soviet Union)\  {\bf 153} (1935).

\bibitem{MCGAHAGAN}
McGahagan,~H.:
{\it An approximation scheme for Schr\"odinger maps}.
 Comm.\ Partial\ Differential\ Equations {\bf 32} (2007), 375--400.


\bibitem{NSU1}
Nahmod,~A., Stefanov,~A., Uhlenbeck,~K.:
{\it On Schr\"odinger maps}.
 Comm.\ Pure Appl.\ Math. {\bf 56} (2003), no. 1, 114--151.

\bibitem{NSU2}
Nahmod,~A., Stefanov,~A., Uhlenbeck,~K.:
{\it Erratum: "On Schr\"odinger maps" [Comm.\ Pure Appl.\ Math. 56  
(2003), no.1, 114--151; MR1929444]}.
 Comm.\ Pure Appl.\ Math.\  {\bf 57} (2004), 833--839.

\bibitem{Nash}
Nash,~J.:
{\it The imbedding problem for Riemannian manifolds}.
 Ann.\ of Math.\ {\bf 63} (1956), 20--63. 

\bibitem{Nishikawa}
Nishikawa,~S.:
{\it Variational Problems in Geometry}.
Translations \ of \ Mathematical \ Monographs \ vol 205 \
American \ Mathematical \ Society, 2002.

\bibitem{NT}
Nishiyama,~T., Tani,~A.:
{\it Initial and initial-boundary value problems for a vortex filament
	with or without axial flow}.
SIAM J.\ Math.\ Anal.\  {\bf 27} (1996), 1015--1023. 

\bibitem{Onodera1}
Onodera,~E.:
{\it A third-order dispersive flow for closed curves into 
K\"ahler manifolds},
J.\ Geom.\ Anal.\ {\bf 18} (2008), in press.

\bibitem{Onodera2}
Onodera,~E.:
{\it Generalized Hasimoto transform of 
one-dimensional dispersive flows 
into compact Riemann surfaces},
SIGMA\ Symmetry\ Integrability\ Geom.\ Methods\ Appl.\ 
 {\bf 4} (2008), 044, 10~pages. 


\bibitem{PWW2}
Pang,~P.~Y.~Y., Wang,~H.~Y., Wang,~Y.~D.:
 {\it Schr\"odinger flow on Hermitian locally symmetric spaces}.
Comm.\ Anal.\ Geom.\ {\bf 10} (2002), 653--681.

\bibitem{SZ}
Shatah,~J., Zeng,~C.:
 {\it Schr\"odinger maps and anti-ferromagnetic chains}.
Comm.\ Math.\ Phys.\ {\bf 262} (2006), 299--315.


\bibitem{SSB}
Sulem,~P.-L., Sulem,~C., Bardos,~C.:
{\it On the continuous limit for a system of classical spins}.
Comm.\ Math.\ Phys.\  {\bf 107} (1986), 431--454.

\bibitem{TN}
Tani,~A., Nishiyama,~T.:
{\it Solvability of equations for motion of a vortex filament 
with or without axial flow}.
Publ.\ Res.\ Inst.\ Math.\ Sci.\  {\bf 33}  (1997), 509--526.

\bibitem{Tarama}
Tarama,~S.:
{\it Remarks on $L^2$-wellposed Cauchy problem for some dispersive equations}.
J.\ Math.\ Kyoto \ Univ.\  {\bf 37}  (1998), 757--765.

\bibitem{ZS}
Zakharov,~V.~E., Shabat,~A.~B.:
{\it Exact theory of two-dimensional self-focusing
and one-dimensional self-modulation of waves in nonlinear media}.
 Soviet Phys.\ JETP.\ 34.\ (1972), 62--69. 
\end{thebibliography}
\end{document}